\documentclass[12pt]{article}

\usepackage{amsmath}
\usepackage{amssymb}
\usepackage{amscd}
\usepackage{epsfig}

 \newcommand{\N}{\mathbb N}
 \newcommand{\qed}{\hfill $\square$}
 \newcommand{\bee}{\begin{equation}}
 \newcommand{\eee}{\end{equation}}
 \newcommand{\Wb}{\mbox {\bf W}}
 \newcommand{\Lb}{\mbox {\boldmath ${\Lambda}$}}
 \newcommand{\Gb}{\mbox {\boldmath ${\Gamma}$}}
 
 \newcommand{\Lbs}{\mbox{\scriptsize\boldmath ${\Lambda}$}}
 
 \newcommand{\Pb}{\mbox {\bf P}}
 \newcommand{\Pbs}{\mbox {\scriptsize{\bf P}}}

 \newcommand{\Ib}{\mbox {\bf I}}
\newcommand{\diam}{\mbox{\rm diam}}
\newcommand{\be}{\begin{eqnarray}}
\newcommand{\ee}{\end{eqnarray}}
\newcommand{\supp}{\mbox{\rm supp}}

\newcommand{\dens}{\mbox{\rm dens}}

\newcommand{\freq}{\mbox{\rm freq}}
\newcommand{\Vol}{\mbox{\rm Vol}}
\newcommand{\es}{\emptyset}

\newcommand{\eps}{{\mbox{$\epsilon$}}}
\newcommand{\e}{{\varepsilon}}

\newcommand{\R}{{\mathbb R}}

\newcommand{\Z}{{\mathbb Z}}
\newcommand{\C}{{\mathbb C}}
\newcommand{\Nat}{{\mathbb N}}

\newcommand{\Ak}{{\mathcal A}}

\newcommand{\Dk}{{\mathcal D}}

\newcommand{\bd}{\partial}

\newcommand{\Pk}{{\mathcal P}}

\newcommand{\Gk}{{\mathcal G}}

\newcommand{\Ok}{{\mathcal O}}
\newcommand{\Sk}{{\mathcal S}}
\newcommand{\Tk}{{\mathcal T}}

\newcommand{\Xt}{X_{\Tk}}

\newcommand{\dist}{\mbox{\rm dist}}

\newcommand{\Lam}{{\Lambda}}

\newcommand{\lam}{\lambda}
\newcommand{\Gam}{\Gamma}
\newcommand{\om}{\omega}
\newcommand{\mmin}{\rm min}
\newcommand{\mmax}{\rm max}

\def\lt{\left}
\def\rt{\right}

 \newtheorem{theorem}{Theorem}[section]
 \newtheorem{lemma}[theorem]{Lemma}
 \newtheorem{prop}[theorem]{Proposition}
 \newtheorem{cor}[theorem]{Corollary}

 \newtheorem{defi}[theorem]{Definition}
 \newtheorem{example}[theorem]{Example}
 \newtheorem{remark}[theorem]{Remark}

\numberwithin{equation}{section}

\begin{document}

\bibliographystyle{unset}

\title{{\sc Consequences of Pure Point Diffraction Spectra for
Multiset Substitution Systems}}
\date \today
\maketitle

\vspace{7mm}
 \centerline {
{\sc Jeong-Yup  Lee $^{\,\rm a}$},
{\sc Robert V.\ Moody $^{\,\rm a}$}\footnote{RVM is grateful for the continuing support of an NSERC Operating Grant in this research.} ,
{\sc Boris Solomyak $^{\,\rm b}$}\footnote{BS acknowledges support from
NSF grants DMS 9800786 and DMS 0099814.}
}

\vspace{7mm}

{\small
\hspace*{4em}

\smallskip
\hspace*{4em}
a:  Dept. of Mathematical and Statistical Sciences, University of Alberta, \\
\hspace*{4em}
\hspace*{2.5em} Edmonton, Alberta T6G 2G1, Canada

\smallskip
\hspace*{4em}
b:  Department of Mathematics, University of Washington,\\
\hspace*{4em}
\hspace*{2.5em} Seattle, WA 98195, USA.}

\abstract There is a growing body of results in the theory
of discrete point sets and tiling systems giving conditions
under which such systems are pure point diffractive. Here we look
at the opposite direction: what can we infer about a discrete
point set or tiling, defined through a primitive substitution system,
given that it is pure point diffractive? Our basic objects are
Delone multisets and tilings,
which are self-replicating under a primitive substitution
system of affine mappings with a common expansive map $Q$.
Our first result gives a partial answer to a question of Lagarias and
Wang: we characterize repetitive substitution Delone multisets that can be
represented
by substitution tilings using a concept of ``legal cluster.''
This allows us to move freely between both types of objects.
Our main result is that  for lattice substitution multiset systems (in arbitrary dimensions) being a regular
model set is not only sufficient for having pure point spectrum---a known fact---but is also necessary.

This completes a circle of equivalences relating pure point dynamical and diffraction spectra,
modular coincidence, and model sets for lattice substitution systems begun by the first two authors
of this paper.

\section{Introduction}
One of the forces behind the development of the theory of
aperiodic order is the unexpected ability of aperiodic structures
to display perfect pure point diffraction. Until the mid-1980's, this type
of phenomenon had been
considered, at least physically, to be strictly
indicative of crystals. Its discovery in new non-crystallographic
metallic solids (now called quasicrystals) and in mathematically aperiodic
structures like the Penrose tilings made it clear just how poor our
understanding of diffraction had been.

Over the past 15 years a great number of diffractive aperiodic tilings and
point sets have been discovered, so that by now there is no shortage of
examples. Still, there remains the fundamental and difficult question:
given a pure point diffraction pattern, what can we say about where it came from?  In this generality we have little hope of answering the question,
partly because diffraction is only an average statistical property of
the structure creating it and partly because it is well known that in any case  the diffraction pattern does not contain enough mathematical information to allow for the solution.

However, if {\it{a priori}} we know something of the structure of the object
producing the pattern then it is possible to make progress. This is the
key behind practical algorithms for solving this inverse problem in the field of crystallography. In the case of aperiodic structures the things are, of course, more difficult. However, if we know that the originating structure is a discrete set of points arising, say, by marking the points of a substitution tiling then already we have a good deal of prior information,
and in such a situation we can draw useful conclusions.

This paper addresses this type of problem. The basic objects of our study
are tilings and certain colored Delone point sets (multisets) in Euclidean space $\R^d$,
arising from primitive substitution systems (see Section 2 for the definitions).
The underlying assumption that we make throughout is that for each $R>0$,
the number of configurations, up to translation, of diameter less than $R$
of our tilings or multisets is finite, a condition called {\em finite local complexity}.
The main problem is to relate the diffractive properties, especially pure point diffraction,
with the geometrical properties of our tilings or multisets
and the substitution systems that create them. In the case of multisets which form a partition of a lattice
in $\R^d$ these geometrical properties assume the particularly explicit form of model sets (cut
and project sets).

A substitution is
a finite collection of affine mappings whose linear part is a common
expansive mapping $Q$ on $\R^d$ (see Definition \ref{def-subst} and
\ref{def-subst-mul}).
Substitutions are familiar in tilings and are the basic method of
construction of many of the most famous aperiodic examples.
The idea is even simpler for multisets.

Lagarias and Wang \cite{lawa} have developed a
theory for substitution Delone multisets that, under certain conditions,
allows one to represent them by substitution tilings in which each tile
exactly corresponds to one point in the multisets (such Delone multisets
are called {\em representable}).
Parenthetically we note that this is considerably more subtle than simply performing the Voronoi construction
or one of its relatives, which generally does not create a tiling inheriting
the underlying substitution.
The ``inflated'' tiles should decompose into unions of
tiles exactly. In fact the Lagarias and Wang tiles may have ``fractal'' boundaries, as in Example 3.10
below.

A question raised in \cite{lawa} was to characterize
representable tilings. The first new result of the paper (Theorem 3.7)
gives a partial answer. We introduce a notion
of {\em legal cluster} and prove that under the assumption of repetitivity,
legality of all clusters is equivalent to representability by a tiling.
This result enables us to pass freely between substitution Delone multisets
and substitution tilings, and to infer certain properties of substitution
Delone multisets from corresponding properties of substitution tilings.
For example, we are able to prove that the dynamical systems of
primitive substitution Delone multisets for which every cluster is legal
are uniquely ergodic (Corollary \ref{cor-erg}).

A point set is pure point diffractive, or has a ``perfect diffraction pattern'',  if the  Fourier transform of its autocorrelation measure is a pure point measure. The reader is referred to \cite{Hof} or \cite{LMS1}
for more on this.  In \cite{LMS1}
we proved that, under the assumption of unique ergodicity, pure point diffractivity implies
that the corresponding dynamical system has pure discrete dynamical
spectrum, the converse of an earlier result by Dworkin \cite{Dwor}.
Thus in this paper we are able to restrict ourselves to the framework of the spectra of dynamical systems,
which we do.

Section 4 contains Theorem \ref{thm-pure} and
Corollary \ref{cor-pure1} which provide necessary and sufficient conditions
for pure pointedness in terms of the asymptotic behaviour of the $Q$-iterates
of almost-periods of the system.  The theorem itself is a generalization from
two to arbitrary  dimension of a result in \cite{soltil}.
At the present stage of our knowledge we could not prove these conditions
directly, without passing to the associated substitution tilings.

Section 5 contains our main results which are an application of the foregoing theory
to the case of {\em lattice} substitution systems.  These are substitution multisets for which
the multiset is a partition of a lattice .
The objective here is to relate pure pointedness to the theory of model
sets (or cut and project sets).
We prove that being a regular model set is not only
sufficient for having pure point spectrum---a known fact---but is also necessary
for certain lattice substitution multiset systems (in arbitrary dimensions).
This completes a circle of equivalences begun in \cite{LM}.

The final result (Theorem \ref{th-main}) links modular coincidence,
the existence of a model set interpretation, the two notions of pure point
spectrum, and the aforementioned condition on almost-periods.
The main parts of this that are new to this
paper are Theorem \ref{th-modul}  and the reworking of the theory of lattice
substitutions begun in \cite{LM} to take account of the more refined
definition of modular coincidence required here.
Theorem \ref{th-main} provides, in particular,
an extension of Dekking's well-known criterion for
pure point spectrum in terms of coincidences,
generalizing it from its original one dimensional
setting in constant-length alphabetic substitutions.
A number of interesting tilings fall into the setting of lattice
substitutions, including the Robinson \cite{Robin}, sphinx, and chair tilings \cite{LM}.

The paper concludes with an Appendix of two subsections.
One establishes the uniform patch frequency property of fixed points of
primitive tiling substitution systems. This is already known in various forms,
but we  offer a proof here since it is hard to find a version in literature
that is fitted to our present needs.
The uniform patch frequency property is equivalent to the
unique ergodicity of the corresponding dynamical system,
which is essential to our arguments on pure point diffractivity.
The second part of the Appendix proves the necessity
in Theorem \ref{thm-pure},
which is a modification of the argument in \cite{soltil}.

The authors are grateful to the referees whose comments helped
to clarify the structure of this paper.

\section{Definitions and notation}

\subsection{Delone multisets}

\noindent
A {\em multiset \footnote{Caution : In \cite{lawa}, which we occasionally
cite below, the word multiset refers to a set with multiplicities.}} or {\em $m$-multiset} in $\R^d$ is a
subset $\Lb = \Lam_1 \times \dots \times \Lam_m
\subset \R^d \times \dots \times \R^d$ \; ($m$ copies)
where $\Lam_i \subset \R^d$. We also write
$\Lb = (\Lam_1, \dots, \Lam_m) = (\Lam_i)_{i\le m}$.
Recall that a Delone set is a relatively dense and uniformly discrete
subset of $\R^d$.
We say that $\Lb=(\Lambda_i)_{i\le m}$ is a {\em Delone multiset} in $\R^d$ if
each $\Lambda_i$ is Delone and $\supp(\Lb):=\bigcup_{i=1}^m \Lambda_i
\subset \R^d$ is Delone.

Although $\Lb$ is a product of sets, it is convenient to think
of it as a set with types or colors, $i$ being the
color of points in $\Lambda_i$.
A {\em cluster} of $\Lb$ is, by definition,
a family $\Pb = (P_i)_{i\le m}$ where $P_i \subset \Lambda_i$ is
finite for all $i\le m$.
Many of the clusters that we consider have the form
$A\cap \Lb := (A\cap \Lambda_i)_{i\le m}$, for a bounded set
 $A\subset \R^d$.
There is a natural
translation $\R^d$-action on the set of Delone multisets and their clusters
in $\R^d$.
The translate of a cluster $\mbox{\bf P}$ by $x \in \R^d$ is
$x + \Pb = (x + P_i)_{i\le m}$.
We say that two clusters $\Pb$ and $\Pb'$ are {\em translationally equivalent}
if $\Pb=x+\Pb'$ for some $x \in \R^d$.
For any two Delone $m$-multisets $\Lb, \Gb$
we define $\Lb \cap \Gb = (\Lam_i \cap \Gam_i)_{i \le m}$ and
$\Lb \triangle \Gb = (\Lam_i \triangle \Gam_i)_{i \le m}$, where $\Lam_i \triangle \Gam_i = (\Lam_i \backslash \Gam_i) \cup
(\Gam_i \backslash \Lam_i)$.
We write $B_R(y)$ for the {\em closed} ball of radius $R$
centered at $y$ and use also $B_R$ for $B_R(0)$.

\begin{defi} \label{def-flc}
A Delone multiset $\Lb$ has {\em finite local complexity
(FLC)} if for every $R>0$ there exists a finite set $Y\subset \supp(\Lb)=
\bigcup_{i=1}^m \Lam_i$  such that
$$
\forall x\in \supp(\Lb),\ \exists\, y\in Y:\
B_R(x) \cap \Lb = (B_R(y) \cap \Lb) + (x-y).
$$
In plain language, for each radius $R > 0$ there are only finitely many
translational classes of clusters whose support lies in some ball of
radius $R$.
\end{defi}

In this paper we will usually assume that our Delone multisets
have FLC.

\medskip

\begin{defi} \label{def-repetitive}
A Delone multiset $\Lb$ is {\em repetitive} if for any compact set
$K \subset \R^d$, $\{t \in \R^d : \Lb \cap K = (t + \Lb) \cap K\}$ is relatively dense; i.e.
there exists $R = R(K) > 0$ such that every open ball
$B_R(y)$ contains at least one element of $\{t \in \R^d : \Lb \cap K = (t + \Lb) \cap K\}$.
\end{defi}

For a cluster $\Pb$ and a bounded set $A\subset \R^d$ denote
$$
L_{\Pbs}(A) = \sharp\{x\in \R^d:\ x+\Pb \subset A\cap \Lb\},
$$
where $\sharp$ means the cardinality.
In plain language, $L_{\Pbs}(A)$ is the number of translates of $\Pb$
contained in $A$, which is clearly finite.
For a bounded set $F \subset \R^d$ and $r > 0$, let
\[
\begin{array}{l}
F^{+r} := \{x \in \R^d:\,\dist(x,F) \le r\},\\
F^{-r} := \{x \in F:\, \dist(x,\partial F) \ge r\} \supset F \setminus
(\partial F)^{+r}.
\end{array}
\]
A {\em van Hove sequence} for $\R^d$ is a sequence
$\mathcal{F}=\{F_n\}_{n \ge 1}$ of bounded measurable subsets of
$\R^d$ satisfying
\be \label{Hove}
\lim_{n\to\infty} \Vol((\partial F_n)^{+r})/\Vol(F_n) = 0,~
\mbox{for all}~ r>0.
\ee

Throughout this paper we deal with concepts that depend on some sort of averaging sequence for their very definition: densities, frequencies, autocorrelation and diffraction measures. Even when not explicitly mentioned, we will always have in mind that these concepts have been defined in terms of some predetermined van Hove sequence $\{F_n\}_{n \ge 1}$.

\begin{defi} \label{def-ucf}
Let $\{F_n\}_{n \ge 1}$ be a van Hove sequence.
The Delone multiset $\Lb$ has {\em uniform cluster frequencies} (UCF)
(relative to $\{F_n\}_{n \ge 1}$) if for any cluster $\Pb$, there is the limit
$$
\freq(\Pb,\Lb) = \lim_{n\to \infty} \frac{L_{\Pbs}(x+F_n)}{\Vol(F_n)} \ge 0,
$$
uniformly in $x\in \R^d$.
\end{defi}

For any subset $\Lb' \subset \Lb$, we define
$$\dens(\Lb') :=
\lim_{n \to \infty} \frac{\sharp(\Lb' \cap F_n)}{\Vol(F_n)} \,,
$$
if the limit exists.

Let $\Lb$ be a Delone multiset and let $X$ be the collection of all
Delone multisets each of whose
clusters is a translate of a $\Lb$-cluster. We introduce a metric
on Delone multisets in a simple variation of the standard way:
 for Delone multisets $\Lb_1$, $\Lb_2 \in X$,
\be \label{metric-multisets}
d(\Lb_1,\Lb_2) := \min\{\tilde{d}(\Lb_1,\Lb_2), 2^{-1/2}\}\, ,
\ee
where
\be
\tilde{d}(\Lb_1,\Lb_2)
&=&\mbox{inf} \{ \e > 0 : \exists~ x,y \in B_{\e}(0), \nonumber \\ \nonumber
&  & ~~~~~~~~~~ B_{1/{\e}}(0) \cap (-x + \Lb_1) = B_{1/{\e}}(0)
\cap (-y + \Lb_2) \}\,.
\ee
For a proof of this, see \cite{LMS1}.

We define $X_{\Lbs} := \overline{\{-h + \Lb : h \in \R^d \}}$ with the metric
$d$.
In spite of the special role
played by $0$ in the definition of $d$, any other point of $\R^d$ may
be used as a reference point, leading to an equivalent metric and more
 importantly the same topology on $X_{\Lbs}$.

The group $\R^d$ acts on $X_{\Lbs}$ by translations which are obviously
homeomorphisms, and we get a topological dynamical system $(X_{\Lbs},\R^d)$.

\subsection{Tilings}

This subsection briefly reviews the basic definitions of tilings and their associated dynamical
systems.
We begin with a set of types (or colors) $\{1,\ldots,m\}$,
which we fix once and for all.
A {\em tile} in $\R^d$ is defined as a pair $T=(A,i)$ where $A=\supp(T)$
(the support of $T$) is a compact
set in $\R^d$ which is the closure of its interior, and
$i=l(T)\in \{1,\ldots,m\}$
is the type of $T$. We let $g+T = (g+A,i)$ for $g\in \R^d$. We say that
a set $P$ of tiles is a {\em patch} if the number of tiles in $P$ is
finite and the tiles of $P$ have mutually disjoint
interiors (strictly speaking, we have to say ``supports of tiles,'' but this
abuse of language should not lead to confusion). The {\em support of
a patch} is the union of the supports of the tiles that are in it.
Note that the support
of a patch need not be connected. The
{\em diameter of a patch} is the diameter of its support. The
{\em translate of a patch} $P$ by $g\in \R^d$ is $g+P := \{g+T:\ T\in P\}$.
We say that two patches $P_1$ and $P_2$ are {\em translationally equivalent}
if $P_2 = g+P_1$ for some $g\in \R^d$.
A tiling of $\R^d$ is a set $\Tk$ of tiles such that
$\R^d = \cup \{\supp(T) : T \in \Tk\}$ and distinct tiles have disjoint
interiors.
Given a tiling $\Tk$, finite sets of tiles of $\Tk$ are called
$\Tk$-patches.

We define FLC, repetitivity, and uniform patch frequencies (UPF), which is
the analog of UCF, on tilings
in the same way as the corresponding properties on Delone multisets.

We always assume that
\begin{itemize}
\item any two $\Tk$-tiles with the same color are translationally equivalent.\\
(Hence there are finitely many $\Tk$-tiles up to translation.)
\item the tiling $\Tk$ has {\em finite local complexity} (FLC), that is,
for any $R>0$ there are finitely many $\Tk$-patches of diameter less
than $R$ up to translation equivalence.
\end{itemize}

Let $X_{\Tk} = \overline{\{-g+\Tk : g \in \R^d\}}$, where $X_{\Tk}$ carries a
well-known topology, given analogously to (\ref{metric-multisets})
for $X_{\Lbs}$,
relative to which it is compact (equivalent to FLC). We have a natural action of
$\R^d$ on $X_{\Tk}$ which makes it a topological dynamical system.
The set $\{- g+\Tk:g \in \R^d \}$ is the orbit of $\Tk$.
A dynamical system is {\em minimal} if the orbit of every element of
$X_{\Tk}$ is dense in $X_{\Tk}$. The minimality of dynamical system
$(X_{\Tk}, \R^d)$ is equivalent to the repetitivity of $\Tk$
(see \cite{Fur}).

\begin{defi} \label{def-cyl2}
Let $P$ be a patch of $\Tk$ or some translate of $\Tk$, and let
 $W \subset \R^d$ be a Borel set.
Define the cylinder set $X_{P,W}$ by
\[X_{P,W} := \{ \Sk \in X_{\Tk}:\, -g+P \mbox{ is an } \Sk\mbox{-patch for some }
g\in W\}.\]
\end{defi}

\medskip

\section{Substitution Systems}

\subsection{Tile-substitutions}

We say that a linear map $Q : \R^d \rightarrow \R^d$ is {\em expansive}
if there is
a $c > 1$ with
\be
d(Qx,Qy) \geq c \cdot d(x,y) \label{def-expan}
\ee
for all $x,y \in \R^d$
and some metric $d$ on $\R^d$ compatible with the standard topology.
This is equivalent to saying that all the eigenvalues of $Q$
 lie outside the closed unit disk in $\C$.

\begin{defi}\label{def-subst}
Let $\Ak = \{T_1,\ldots,T_m\}$ be a finite set of tiles in $\R^d$
such that $T_i=(A_i,i)$; we will call them {\em prototiles}.
Denote by $\Pk_{\Ak}$ the set of
patches made of tiles each of which is a translate of one of $T_i$'s.
We say that $\omega: \Ak \to \Pk_{\Ak}$ is a {\em tile-substitution} (or simply
{\em substitution}) with
expansive map $Q$ if there exist finite sets $\Dk_{ij}\subset \R^d$ for
$i,j \le m$, such that
\begin{equation}
\om(T_j)=
\{u+T_i:\ u\in \Dk_{ij},\ i=1,\ldots,m\} \ \ \  \mbox{for} \  j\le m,
\label{subdiv}
\end{equation}
with
$$
Q(A_j) = \bigcup_{i=1}^m (\Dk_{ij}+A_i).
$$
Here all sets in the right-hand side must have disjoint interiors;
it is possible for some of the $\Dk_{ij}$ to be empty.
\end{defi}

The substitution (\ref{subdiv}) is extended to all translates of prototiles by
$\om(x+T_j)= Q x + \om(T_j)$, and to patches and tilings by
$\om(P)=\cup\{\om(T):\ T\in P\}$.
The substitution $\om$ can be iterated, producing larger and larger patches
$\om^k(T_j)$. To the substitution $\om$ we associate its $m \times m$
substitution matrix $S$, with $S_{ij}:=\sharp (\Dk_{ij})$.
The substitution $\om$ is called {\em primitive}
if the substitution matrix $S$ is primitive, i.e.\ there is an $l > 0$ for
which $S^l$ has no zero entries.

\begin{defi} \label{def-selfaf}
A patch will be called {\em legal} if it is a translate
of a subpatch of $\om^k(T_i)$ for some $i\le m$ and $k\ge 1$.
A tiling $\Tk$ with FLC is said to be {\em self-affine}
with the prototile set $\Ak$, expansive map $Q$, and primitive
substitution $\om$, if every $\Tk$-patch is  legal.
\end{defi}

The set of self-affine tilings associated with $(\Ak,\om)$ will be
denoted by $X_{\Ak,\om}$.
A tiling $\Tk$ is called a {\em fixed point}
of the substitution $\om$ if $\om(\Tk)=\Tk$. It turns out that such a
tiling need not be repetitive, even though the substitution is primitive,
see \cite{solcorr}.
In that case one can show that
$X_{\Tk} \supset X_{\Ak,\om}$ but this inclusion is proper.
It is well-known (and easy to see) that one can always find a {\em periodic
point} for $\om$ in the space $X_{\Ak,\om}$ i.e. there is $\Tk \in
X_{\Ak,\om}$ such that $ \om^N (\Tk) = \Tk$ for some $N \ge 1$.
In this case we can always use $\omega^N$ in place of $\omega$ to obtain a tiling which is a fixed point.

\begin{lemma} \label{rep-legal}
Let $\Tk$ be a fixed point of a primitive substitution $\om$
with expansive map $Q$ and prototiles $\Ak$.
Then $\Tk$ is repetitive if and only if every $\Tk$-patch is legal, i.e.
$\Tk \in X_{\Ak,\om}$.
\end{lemma}

\noindent
{\em Proof.}
Suppose $\Tk$ is repetitive. Then for any patch $P$ of $\Tk$ there exists
$R > 0$ such that every open ball $B_R(y)$ contains a patch translationally
equivalent to $P$. Since each tile of $\Tk$ has non-empty interior and $\om$
 is with expansive map $Q$, there is $M \ge 1$ such that for any $i \le m$,
$Q^M(A_i)$ contains $B_R(y)$ for some $y \in \R^d$. This means that
$\om^M(T_i)$ contains a patch translationally equivalent to $P$ for any
$i \le m$. Thus every $\Tk$-patch is legal.

Conversely, suppose every $\Tk$-patch is legal. Then for every $\Tk$-patch $P$
there is $K \ge 1$ such that for any $i \le m$, $P$ is a translate of
a subpatch of $\om^K(T_i)$ by the primitivity of $\om$.
Choose $r > {\mmax}\{{\diam(T_i)}:\  i \le m \}$.
Every open ball $B_r(y)$ contains at least one tile.
So every open ball $B_{\|Q\|^K r}(y)$, where $\|Q\|$ is the operator norm,
contains at least one supertile
$\om^K(T)$, which contains a translate of $P$.
Therefore $\Tk$ is repetitive. \qed

\subsection{From substitution Delone multisets to tilings.} \label{Delone-tiling}

We now link the theory of multisets and tilings through the notion of
representable Delone multisets.
The concept was introduced by Lagarias and Wang
\cite{lawa}, under the name of self-replicating Delone set families.

\begin{defi} \label{def-subst-mul}
$\Lb = (\Lam_i)_{i\le m}$ is called a {\em
substitution Delone multiset} if $\Lb$ is a Delone multiset and
there exist an expansive map
$Q:\, \R^d\to \R^d$ and finite sets $\Dk_{ij}$ for $i,j\le m$ such that
\be \label{eq-sub}
\Lambda_i = \bigcup_{j=1}^m (Q \Lambda_j + \Dk_{ij}),\ \ \ i \le m,
\ee
where the unions on the right-hand side are disjoint.
\end{defi}

We say that the
substitution Delone multiset is {\em primitive} if the corresponding
substitution matrix $S$, with $S_{ij}= \sharp (\Dk_{ij})$, is primitive.

Let $Y$ be a nonempty set in $\R^d$.
For $m \in \Z_{+}$, an $m \times m$ {\em matrix function system} (MFS)
 on $Y$ is an $m \times m$ matrix $\Phi = (\Phi_{ij})$, where each $\Phi_{ij}$
is a finite set (possibly empty) of mappings $Y$ to $Y$.

Any MFS $\Phi$ on $Y$ induces a mapping on Delone multiset
$\Lb = (\Lambda_i)_{i \le m}$, where $\Lambda_i \subset Y$, $i \le m$, by
\begin{eqnarray}
{\Phi \left[ \begin{array}{c}
                     \Lambda_{1} \\ \vdots \\ \Lambda_{m}
             \end{array}
      \right]}&=&{\left[ \begin{array}{c}
                         \bigcup_{j \le m} \bigcup_{f \in \Phi_{1j}} f(\Lambda_{j})\\ \vdots \\
 \bigcup_{j \le m} \bigcup_{f \in \Phi_{mj}} f(\Lambda_{j})
                         \end{array}
                  \right]_{\,,}}
\end{eqnarray}
which we call the {\em substitution determined by $\Phi$}.
We often write $\Phi_{ij} (\Gamma_{j})$ for
$\bigcup_{f \in \Phi_{ij}} f(\Gamma_{j})$, $\Phi(\Gamma_j)$
for $(\Phi_{ij} (\Gamma_{j}))_{i \le m}$, and
$\Phi(\Gb)$ for $(\bigcup_{j \le m}\Phi_{ij} (\Gamma_{j}))_{i \le m}$
for any subset $\Gb = (\Gamma_j)_{j \le m} \subset \Lb$.
In particular, we often write $\Phi_{ij}(x)$ for $\Phi_{ij}(\{x\})$,
where $x \in \Lam_j$ .
We associate to $\Phi$ its substitution matrix $S(\Phi)$, with
$(S(\Phi))_{ij}= \sharp (\Phi_{ij})$.

Let $\Phi, \Psi$ be $m \times m$ MFS's on $Y$.
Then we can compose them :
\be
\Psi \circ \Phi = ((\Psi \circ \Phi)_{ij})\,,
\ee
where  $(\Psi \circ \Phi)_{ij} = \bigcup_{k=1}^{m} \Psi_{ik} \circ
\Phi_{kj}~$and$~\Psi_{ik} \circ \Phi_{kj} :=
\left\{ \begin{array}{l} \{~g \circ f~ :~ g \in \Psi_{ik}, f \in \Phi_{kj}~\} \\
 \,\emptyset ~~~~~~ \mbox{if}~ \Psi_{ik} = \emptyset~ \mbox{or}~ \Phi_{kj} = \emptyset
\,.                      \end{array}
  \right.$
Evidently, $S(\Psi \circ \Phi) \leq S(\Psi)\, S(\Phi)$.

For any given substitution Delone multiset $\Lb = (\Lambda_i)_{i \le m}$,
we can always find the corresponding MFS $\Phi$ such that $\Phi(\Lb)=\Lb$.
Indeed, by Definition \ref{def-subst-mul} we can take
$\Phi_{ij} = \{ f : \, f : x \mapsto Qx + a,\,a \in \Dk_{ij}\}$.
So $\Phi(\Lam_j) = (Q \Lam_j + \Dk_{ij})_{i \le m}$, $j \le m$.

\begin{lemma}{\em \cite[Th.2.3]{lawa}}
If $\Lb$ is a primitive substitution Delone multiset with expansive map $Q$,
then the PF eigenvalue of its substitution matrix $S$ equals $|\det(Q)|$.
\end{lemma}
\medskip

For each primitive substitution Delone multiset there is an
{\em adjoint system} of equations
\be \label{eq-til}
Q A_j = \bigcup_{i=1}^m (\Dk_{ij} + A_i),\ \ \ j \le m.
\ee
From Hutchinson's Theory (or rather, its generalization to the
``graph-directed''
setting), it follows that (\ref{eq-til}) always has a unique solution
for which $(A_i)_{i \le m}$ is
a family of non-empty compact sets of $\R^d$
(see for example \cite{BM1}, Prop.1.3).
It is proved in \cite[Th.2.4 and Th.5.5]{lawa} that if $\Lb$ is a primitive
substitution Delone multiset, then all the sets $A_i$ from (\ref{eq-til})
have non-empty interior and, moreover, each $A_i$ is the closure of
its interior.

\begin{defi} A Delone multiset $\Lb = (\Lam_i)_{i \le m}$ is called
{\em representable} (by tiles) if there exist tiles $T_i = (A_i,i), i\le m, $
so that
\be\label{eq-1}
\{x + T_i :\ x\in \Lambda_i,\ i \le m\} \ \ \ \mbox{is a tiling of}\ \
\R^d,
\ee
that is, $\R^d = \bigcup_{i\le m} \bigcup_{x\in \Lambda_i} (x + A_i)$, and the
sets in this union have disjoint interiors.
In the case that $\Lb$ is a primitive substitution Delone multiset we will
understand the term representable to mean relative to the tiles
$T_i = (A_i,i), i\le m$, arising from the solution to the adjoint system (\ref{eq-til}).
\end{defi}

A cluster will be called {\em legal} if it is a translate of a subcluster of
$\Phi^k(x_j)$ for some $x_j \in \Lam_j$, $j \le m$ and some $k \in \Z_+$.

\medskip

Lagarias and Wang have a condition, namely existence of a fundamental cycle
of period $1$, which ensures representability. Legality generalizes this and
in fact our Theorem \ref{legal-rep} is based on \cite[Th.7.1]{lawa}.

In the same paper \cite[Lemma 3.2]{lawa} it is shown that if $\Lb$ is a substitution Delone multiset,
then there is a finite multiset (cluster) $\Pb \subset \Lb$ for which
$\Phi^{n-1}(\Pb) \subset \Phi^n(\Pb)$ for $n \ge 1$ and
$\Lb = \lim_{n \to \infty} \Phi^n (\Pb)$. We call such a multiset $\Pb$
a {\em generating multiset}.

\begin{theorem}\label{legal-rep}
Let $\Lb$ be a primitive substitution Delone multiset such that every
$\Lb$-cluster is legal.
Then $\Lb$ is representable.
\end{theorem}

\noindent {\em Proof.}
Let $\Phi$ be the MFS satisfying $\Phi(\Lb)=\Lb$ and suppose
$\Ak = \{T_1, \dots, T_m\}$ arises from
the solution to the adjoint system (\ref{eq-til}).
So for any $M \in \Z_+$,
\be \label{adj-equation}
 Q^M(A_j) = \bigcup_{i=1}^m ((\Dk^M)_{ij} + A_i), \ \ \ j \le m,
\ee
where $A_i = \supp(T_i)$ and
\[(\Dk^M)_{ij} = \bigcup_{k_1,k_2,\dots,k_{(M-1)} \le m}
(\Dk_{ik_1} + Q \Dk_{k_1 k_2} + \cdots + Q^{M-1} \Dk_{k_{(M-1)} j}).
\]
On the other hand, for any $M \in \Z_+$ and $i \le m$,
\be \label{substi-formula}
(\Phi^M)_{ij} (x_j) = Q^M x_j +(\Dk^M)_{ij}, \ \ \mbox{for any} \
x_j \in \Lam_j, j \le m.
\ee
So putting (\ref{adj-equation}) and (\ref{substi-formula}) together, we get
\be \label{adj-equation-substi}
Q^M(x_j &+& A_j)  =  \bigcup_{i=1}^m ( Q^M x_j +(\Dk^M)_{ij} + A_i)
 \nonumber \\
                & = & \bigcup_{i=1}^m ((\Phi^M)_{ij} (x_j) + A_i)
 \nonumber \\
               & = & \bigcup_{i=1}^m \bigcup_{y \in (\Phi^M)_{ij} (x_j)}
(y + A_i),  \ \ \mbox{for any} \ x_j \in \Lam_j, \  j \le m.
\ee
From \cite[Th.2.4 and Th.5.5]{lawa}, $\mu(A_j) > 0$ and $A_j$ is the
closure of its interior for any $j \le m$.
Let $\tilde\mu := (\mu(A_1), \dots, \mu(A_m)) > 0$.
Taking measures in (\ref{adj-equation-substi}),
\be \label{inequ-measure}
 |\det Q|^M \tilde\mu \le \tilde\mu S(\Phi^M) \le \tilde\mu S(\Phi)^M,
\ee
where $S(\Phi)$ is the substitution matrix of $\Phi$.
From \cite[Th.2.3]{lawa}, we know that
$|\det Q| = $ PF eigenvalue of $S(\Phi)$.
So we can derive
\[ |\det Q|^M \tilde\mu = \tilde\mu S(\Phi)^M
\]
from (\ref{inequ-measure}), see \cite[Lemma 1]{LM}.
Thus for any $x_j \in \Lam_j, j \le m$,
\be \label{measure-adj-equation}
\mu(Q^M(x_j + A_j)) &=& \mu \left(\bigcup_{i=1}^m ((\Phi^M)_{ij} (x_j) + A_i)\right)
\nonumber \\
&=& \sum_{i=1}^m (S(\Phi^M))_{ij} \mu(A_i),
\ee
and this shows that the unions on the right-hand side of
(\ref{adj-equation-substi}) are measure-wise disjoint.
So we get a (tile)-substitution $\om : \Ak \to \Pk_{\Ak}$ associated
with $\Ak$ and $\Lb$.

Next, let $\Pb$ be a generating multiset. From the assumption that every
cluster in $\Lb$ is legal, there is $a_l \in \Lam_l$ for some $l \le m$
such that $z + \Pb = (z + P_j)_{j \le m} \subset
((\Phi^K)_{jl}(a_l))_{j \le m}$
 for some $K \in \Z_+$ and some $z \in \R^d$.
So
\be \label{legal-patch}
 \{ x + T_j : x \in (\Phi^K)_{jl}(a_l), j \le m\} \supset
 \{ x + T_j : x \in (z + P_j), j \le m\}.
\ee
By (\ref{measure-adj-equation}), all the tiles in the left-hand side of
(\ref{legal-patch}) are measure-wise disjoint.
So then are the tiles in the right-hand side of (\ref{legal-patch}).
This implies that all the tiles in the set
$\{p_j + T_j : p_j \in P_j, j \le m\}$ are measure-wise disjoint.
Thus for any $n \in \Z_+$, all the super tiles in the set
$\{\om^n(p_j + T_j) : p_j \in P_j, j \le m\}$ are measure-wise disjoint.
Noting that
$$\om^n(p_j + T_j) = \{ x+T_i : x \in (\Phi^n)_{ij}(p_j), i \le m\} \ \
\mbox{for each} \ p_j \in P_j,\ j \le m,
$$
we get that
$$\Phi^n(\Pb) + \Ak := \{ x+T_i : x \in (\Phi^n)_{ij}(p_j), p_j \in P_j,
i,j \le m \}
$$
also consists of tiles which are measure-wise disjoint.
Since $\Phi^{n-1}(\Pb) \subset \Phi^n(\Pb)$ for $n \ge 1$ and
$\Lb = \lim_{n \to \infty} \Phi^n (\Pb)$,
$\Lb + \Ak := \{x_j + T_j : x_j \in \Lam_j, j \le m\}$ consists of tiles
which are measure-wise disjoint. Thus distinct tiles in $\Lb + \Ak$ have
disjoint interiors.

Now we need prove that $\Lb + \Ak$ is a tiling.
Let us define $\supp(\Lb + \Ak):= \cup \{x_j + A_j :
x_j \in \Lam_j, j \le m\}$.
Then
\be \label{invariant-underQ}
Q(\supp(\Lb + \Ak)) &=& \cup \{Qx_j + QA_j : x_j \in \Lam_j, j \le m\}
\nonumber \\
&=& \cup \{Qx_j + \Dk_{ij} + A_i : x_j \in \Lam_j, i,j \le m\}
\nonumber \\
&=& \cup \{\Phi_{ij}(x_j) + A_i : x_j \in \Lam_j, i,j \le m\}
\nonumber \\
&=&  \cup \{x_i + A_i : x_i \in \Lam_i, i \le m\}
\nonumber \\
&=& \supp(\Lb + \Ak).
\ee
Suppose $\R^d \backslash \supp(\Lb + \Ak) \neq \emptyset$.
Then there is $z \in \R^d \backslash \supp(\Lb + \Ak)$ and a ball
$B_r(z)$ with radius
$r$ centered at $z$ such that  $B_r(z) \subset \R^d \backslash \supp(\Lb + \Ak)$,
since $\supp(\Lb + \Ak)$ is a closed set.
So for any $N \in \Z_+$,
\be
Q^N(B_r(z)) \cap \supp(\Lb + \Ak) &=& Q^N(B_r(z)) \cap Q^N(\supp(\Lb + \Ak))
\ \ \ \mbox{by (\ref{invariant-underQ})} \nonumber \\
&=&  Q^N(B_r(z) \cap \supp(\Lb + \Ak)) = \emptyset \nonumber.
\ee
But this is a contradiction, since $\Lb$ is Delone.
Therefore $\Lb + \Ak$ is a tiling and so $\Lb$ is representable.
\qed

\begin{cor} \label{leg-rep}
Let $\Lb$ be a repetitive primitive substitution Delone multiset.
Then every $\Lb$-cluster is legal if and only if $\Lb$ is representable.
\end{cor}

\noindent {\em Proof.}
We only need to prove the sufficiency direction.
Suppose $\Lb$ is representable. Then we have a tiling
$\Lb + \Ak = \{x_i + T_i : x_i \in \Lam_i, i \le m\}$ with a
unique solution $\Ak = \{T_1, \dots, T_m\}$ of the adjoint system of
equations such that $\om(\Lb + \Ak) = \Lb + \Ak$, where $\om : \Ak \to
\Pk_{\Ak}$ is a (tile)-substitution. So $\Lb + \Ak$ is a fixed point of
a primitive substitution $\om$ with expansive map $Q$.
By Lemma \ref{rep-legal}, every $(\Lb + \Ak)$-patch is legal, since $\Lb + \Ak$ is repetitive.
Recall that for any $M \in \Z_+$,
$$\om^M(x_j + T_j) =
\{ x + T_i : x \in (\Phi^M)_{ij}(x_j), i \le m\} \ \ \mbox{for any} \
x_j \in \Lam_j,\ j \le m,
$$
where $\Phi$ is the MFS satisfying $\Phi(\Lb) = \Lb$.
So the legality in the tiling $\Lb + \Ak$ shows the legality in $\Lb$.
Therefore every $\Lb$-cluster is legal.
\qed

\medskip

\begin{remark}
{\em Note that, in order to check that every $\Lb$-cluster is legal,
we only need to see if some cluster that contains a finite generating
multiset for $\Lb$ is legal.}
\end{remark}

\medskip

Any tiling $\Tk$ can be converted into a Delone multiset by simply choosing
a point $x_{(A,i)}$ for each tile $(A,i)$ so that the chosen
points for tiles of the same type are in the same relative position in the tile:
$x_{(g+A,i)}= g + x_{(A,i)}$.
We define $\Lambda_i := \{x_{(A,i)}:\ (A,i) \in \Tk \}$ and
$\Lb := (\Lambda_i)_{i \le m}$. Clearly $\Tk$ can be reconstructed from $\Lb$
given the information about how the points lie in their respective tiles.
This bijection establishes a topological conjugacy of
$(X_{\Lbs},\R^d)$ and $(X_{\Tk},\R^d)$. Concepts and theorems can clearly be
interpreted in either language (FLC, UCF,  unique ergodicity,
pure point dynamical spectrum, etc.).

In the case that $\Lb$ is a representable primitive substitution
Delone multiset, we make use of the following bijection.
We consider $T_i = (A_i,i)$, $i\le m$, as prototiles, where
$A_i$'s are defined by (\ref{eq-til}).
Let $\Tk = \Tk(\Lb)$ be the tiling in (\ref{eq-1}), with the colors added,
that is, $\Tk = \{x_i + T_i :\ x_i \in \Lambda_i,\ i\le m\}$,
and let $\Ak = \{T_1,\ldots,T_m\}$.
By (\ref{eq-til}) and the definition of representable primitive substitution
Delone multiset, we have a tile-substitution $\om:\ \Ak\to \Pk_{\Ak}$.

\begin{lemma} \label{lem-conj} The dynamical systems $(X_{\Lbs},\R^d)$ and
$(X_{\Tk},\R^d)$ are topologically conjugate.
\end{lemma}

\noindent {\em Proof.} For any Delone multiset $\Gb = (\Gam_i)_{i \le m}
\in X_{\Lbs}$ we let $\phi(\Gb) = \{x_i + T_i:\ x_i \in \Gam_i,\ i\le m\}$.
It is straightforward to check that $\phi(\Gb) \in X_{\Tk}$ and $\phi$ is a
homeomorphism commuting with the translation action.
\qed

\begin{example}
(A substitution Delone multiset in $\R^2$ with gasket tiles.)\\
{\em  Consider the substitution on $\R^2$ with the following MFS $\Phi$;
\begin{eqnarray*}
 \Phi = \left( \begin{array}{cccc}
              \{f_1, f_4\} & \{f_1\} & \{f_1\} & \emptyset \\
              \{f_2\} & \{f_2, f_3\} & \emptyset &
                         \{f_2\}\\
              \{f_3\} & \emptyset & \{f_2, f_3\} &
                         \{f_3\} \\
              \emptyset & \{f_4\} & \{f_4\} &
                          \{f_1, f_4\}
              \end{array}  \right),
\end{eqnarray*}
where $f_1(x)= 2x$, $f_2(x) = 2x+(1,0)$, $f_3(x) = 2x+(0,1)$, and
$f_4(x) = 2x+(-1,-1)$.

The Delone multiset $\Lb = (\Lam_1, \Lam_2, \Lam_3, \Lam_4)$ generated
from ($\{(0,0),(1,1)\}$, $\{(0,-1)\}$, $\{(-1,0)\}$, $\emptyset$) is
fixed under $\Phi$. We observe that $\bigcup_{i \le 4} \Lam_i = \Z^2$.  The generating multiset
($\{(0,0),(1,1)\}$, $\{(0,-1)\}$, $\{(-1,0)\}$, $\emptyset$) is legal, indicating representability. The solution from the adjoint system consists of four copies of the gasket tile($T$) \cite{Vin}.
Thus $\{x + T : x \in \Lam_i, i \le 4\}$ is a tiling of $\R^2$ by gaskets
pinned down on the standard lattice $\Z^2$.  }
\begin{figure}[ht]
\centerline{\epsfysize=30mm \epsfbox{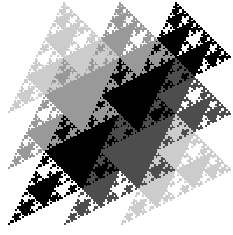}}
\caption{Gasket Tiling}  \label{gasket:1}
\end{figure}
\end{example}

\begin{example} \label{ex-nonlegal-rep}
(A substitution Delone multiset for which not every cluster is legal.)\\
{\em Consider the MFS $\Phi$
\begin{eqnarray*}
\Phi = \left( \begin{array}{cc}
              \{5x, 5x+2, 5x+4\} & \{5x+1, 5x+3\} \\
              \{5x+1, 5x+3\} & \{5x, 5x+2, 5x+4\}
              \end{array}  \right)
\end{eqnarray*}
which generates the bi-infinite sequence shown below with the $a$ and $b$
point sets starting from the generating set $(\{0\}, \{-\frac{1}{2}, -1\})$.
This leads to the
Delone multiset $\Lb = (\Lam_a, \Lam_b)$ which is fixed under $\Phi$.

\begin{picture}(100,50)(0,15)
 \put(5,40){\line(1,0){350}}
 \put(20,40){\multiput(20,0)(20,0){15}{\circle*{.20}}}
 \put(15,30){$\cdots$}
 \put(40,30){$b$}
 \put(60,30){$b$}
 \put(80,30){$a$}
 \put(100,30){$a$}
 \put(120,30){$b$}
 \put(140,30){$b$}
 \put(160,30){$a$}
 \put(180,30){$a$}
 \put(200,30){$b$}
 \put(220,30){$b$}
 \put(240,30){$a$}
 \put(260,30){$a$}
 \put(280,30){$b$}
 \put(300,30){$b$}
 \put(320,30){$a$}
 \put(340,30){$\cdots$}
\put(15,45){$\cdots$}
 \put(40,46){-3}
 \put(60,46){-$\frac{5}{2}$}
 \put(80,46){-2}
 \put(100,46){-$\frac{3}{2}$}
 \put(120,46){-1}
 \put(140,46){-$\frac{1}{2}$}
 \put(160,46){0}
 \put(180,46){$\frac{1}{2}$}
 \put(200,46){1}
 \put(220,46){$\frac{3}{2}$}
 \put(240,46){2}
 \put(260,46){$\frac{5}{2}$}
 \put(280,46){3}
 \put(300,46){$\frac{7}{2}$}
 \put(320,46){4}
 \put(340,45){$\cdots$}
\end{picture}

Note that the generating set $(\{0\}, \{-\frac{1}{2}, -1\})$ is not legal,
even though $\Lb$ is periodic and so repetitive.
Thus $\Lb$ is not representable by tiles arising from the solution to the
adjoint system. In fact the solution to the adjoint system is
$\Ak = \{[0,1],[0,1]\}$ so that $\Lb + \Ak$ double tiles
the line.}
\end{example}

\medskip

\section{Pure point spectrum} \label{section-pure-point}

\subsection{Pure-pointedness for tilings}

This subsection is largely based on \cite{soltil}. The main new feature is
that here the dimension $d$ is arbitrary, while \cite{soltil} was mostly
focused on the case $d=2$.

\smallskip

\medskip
Recall that a topological dynamical system is uniquely ergodic if there is a
unique invariant Borel probability measure.

\begin{theorem} \label{thm-erg}
If $\Tk$ is a fixed point of a primitive substitution with FLC
\footnote{Recently Ludwig Danzer \cite{Danzer} has given an example of a
primitive substitution tiling which does not satisfy FLC.},
then $\Tk$ has uniform patch frequencies, and
the tiling dynamical system $(X_{\Tk}, \R^d)$ is uniquely ergodic.
\end{theorem}

This is already known in various forms
but we offer a proof in the Appendix
since it is hard to find a version in literature
that is fitted to our present needs.
More precisely, in Lemma \ref{lemma-freq} we establish the
existence of uniform patch frequencies (UPF). The fact that
UPF implies unique ergodicity is proved in
\cite[Theorem 2.6 and Theorem 2.7]{LMS1} in the setting of arbitrary
 Delone multisets with FLC.

\medskip

The theorem shows that the dynamical system $(X_{\Tk}, \R^d)$ is ergodic
with respect to the unique invariant probability measure on $X_{\Tk}$, which we
will denote by $\mu$.
We will also need the following result, proved in \cite[Corollary 2.8]{LMS1}.
Let $\eta(\Tk) > 0$ be chosen such that every tile support contains a ball of
diameter $\eta(\Tk)$.

\begin{cor} \label{cor-meas2}
Let $\Tk$ be a fixed point of a primitive substitution with FLC.
Then for any $\Tk$-patch $P$ and
any Borel set $V$ with $\diam(V) <\eta(\Tk)$, we have
$$
\mu(X_{P,V}) = \Vol(V) \cdot\freq(P,\Tk).
$$
\end{cor}

We consider the associated
group of unitary operators $\{U_g\}_{g\in \R^d}$ on $L^2(\Xt,\mu):$
\[U_g f(\Sk) = f(-g+\Sk).\]
A vector $\alpha =(\alpha_1,\ldots,\alpha_d) \in \R^d$ is said to be
an eigenvalue for the $\R^d$-action if there exists an eigenfunction
$f\in L^2(\Xt,\mu),$ that is, $\ f\not\equiv 0$ and
\[U_g f(\Sk) = e^{2 \pi i g \cdot \alpha} f(\Sk),\ \ \ \mbox{for all}
\ \ g\in \R^d.\]

The dynamical system $(\Xt,\mu,\R^d)$ is said to have {\em pure discrete}
(or pure point) {\em spectrum} if the linear span of the
eigenfunctions is dense in $L^2(\Xt)$.

\medskip

Let $\Xi(\Tk)$ be the set of
translation vectors between $\Tk$-tiles of the same type:
\be \label{def-xi}
\Xi(\Tk)= \{x\in \R^d:\ \exists \,T,T' \in \Tk, \ T'=x+T\}.
\ee
Since $\Tk$ has the inflation symmetry with the expansive map $Q$,
we have that
$Q\Xi(\Tk) \subset \Xi(\Tk)$. Note also that $\Xi(\Tk) = -\Xi(\Tk)$.
If $\Tk=\Tk(\Lb)$ is a tiling for a representable Delone multiset $\Lb$,
then $\Xi(\Tk) = \bigcup_{i=1}^m (\Lambda_i-\Lambda_i)$, see Section
\ref{Delone-tiling}.
\medskip
The following is proved in \cite[\S 4]{soltil}.

\begin{theorem} \label{thm-eigen} Suppose that $\Tk$ is a repetitive fixed
point of a primitive substitution with expansive map $Q$ and $\Tk$ has FLC.
If $\alpha \in \R^d$ is an eigenvalue for
$(\Xt,\mu,\R^d)$, then for any $x\in \Xi(\Tk)$ we have
$$
\lim_{n\to\infty}
e^{2\pi i (Q^n x) \cdot \alpha}=1.
$$
\end{theorem}

\medskip

This theorem  yields necessary conditions on the expansive map
$Q$ for the dynamical system to have non-trivial eigenfunctions. The
simplest is that if $Q$ is a diagonal matrix with diagonal entries
$\lam>1$, then $\lam$ has to be a Pisot number. This follows from
the algebraicity of $\lam$ (easy to see) and the classical theorem of Pisot.
Other conditions can be found in \cite{soltil}.

\medskip

For $x\in \Xi(\Tk)$ consider the infinite subset
$$
D_x :=  \Tk \cap (x+\Tk).
$$
It is non-empty by (\ref{def-xi}), and $\supp(D_x)$ is
relatively dense by repetitivity. Observe that $D_x$ has a well-defined
density given by
\be \label{eq-dens}
\dens(D_x) & = & \lim_{n \to \infty}\frac{\Vol(D_x \cap F_n)}{\Vol(F_n)}
\nonumber \\
& = &
\sum_{i=1}^m \freq(T_i\cup (x+T_i),\Tk)\cdot \Vol(A_i) \, > \,0
\ee
where $\{F_n\}_{n \ge 1}$ is a van Hove sequence, $T_i$'s are representatives of all tile types
and $A_i$'s are their supports. For this reason we may call the elements of $\Xi(\Tk)$ the {\em almost-periods} of $\Tk$. Of course an almost-period only
really ``looks like'' a period if the corresponding density $\dens(D_x)$ is
close to $1$.

Below we consider the cylinder set $X_{\{T\},V}$ for a $\Tk$-tile $T$ and
a Borel set $V$, which we denote $X_{T,V}$ to simplify the notation.

\begin{lemma} \label{cyl-int-patch}
Suppose $\Xi(\Tk)$ is uniformly discrete.
Then there is $r > 0$ and $n_0 > 0$ such that for any Borel set $V$ with
$\diam(V) < r$, all $x \in \Xi(\Tk)$, and every $\Tk$-tile $T$,
\be \label{eq-Meyer}
X_{T,V} \cap X_{Q^nx+T,V} =  X_{T\cup (Q^nx+T),V}\,, \ \ \mbox{for all}\
n \ge n_0.
\ee
\end{lemma}
\medskip

\noindent {\em Proof.} Note that for any $x \in \Xi(\Tk)$ there is $n_0 > 0$
such that $T \cup (Q^n x + T)$ is a $\Tk$-patch for all $n \ge n_0$, from
the primitivity.
We only need to check ``$\subset$''
since the other inclusion is obvious in all cases. A tiling
$\Sk\in \Xt$ is in the left-hand side of (\ref{eq-Meyer}) if and only if there
are two, possibly distinct vectors $v_1,v_2\in V$ such that
$-v_1+T\in \Sk$ and $-v_2+Q^n x+T \in \Sk$. Hence
$v_1-v_2+Q^nx \in \Xi(\Sk) = \Xi(\Tk)$. But $Q^n x\in \Xi(\Tk)$ since
$Q\Xi(\Tk) \subset \Xi(\Tk)$. Since $\Xi(\Tk)$ is
uniformly discrete, $v_1-v_2=0$ if
$\diam(V)$ is sufficiently small. Then $T\cup (Q^n x+T) \subset v_1+\Sk$, and
hence $\Sk$ is in the right-hand side of (\ref{eq-Meyer}). The lemma is
proved. \qed

\begin{prop} \label{prop-pure} Suppose that $\Tk$ is a repetitive fixed
point of a primitive substitution with expansive map $Q$
such that $\Xi(\Tk)$ is uniformly discrete and $\Tk$ has FLC.
 If $(\Xt,\mu,\R^d)$ has pure discrete spectrum, then
\be \label{eq-pure1}
\lim_{n\to\infty} \dens(D_{Q^nx})=1,\ \ \ \mbox{for all}\ \  x\in \Xi(\Tk).
\ee
\end{prop}
\medskip

\noindent {\em Proof.}
Fix $x\in \Xi(\Tk)$.
By Theorem~\ref{thm-eigen}, for every eigenvalue $\alpha \in \R^d$ we have
$e^{2\pi i (Q^n x) \cdot \alpha} \to 1.$
This implies
\be \label{eq-conv}
({U}_{Q^nx} -I) f_\alpha \to 0,
\ee
in the norm of $L^2(\Xt,\mu)$, for the corresponding eigenfunction $f_\alpha$.
Since $\|{U}_{Q^nx} -I\| \le 2,$ the sequence of operators
$\{{U}_{Q^nx}-I\}_{n\ge 0}$ is uniformly bounded, so by (\ref{eq-conv}) we have
${U}_{Q^nx} f\to f$ for any $f$ in the closed linear span of the eigenfunctions.
Let $T$ be a $\Tk$-tile and let $V$ be a Borel set satisfying
(\ref{eq-Meyer}).  Denote by $f$
the characteristic function of the cylinder set $X_{T,V}$,
which is in the closed linear span of the eigenfunctions by assumption.
We can write
\be \|{U}_{Q^nx}f-f\|_2^2 & = & \int_{\Xt} |f(-Q^nx+\Sk)-f(\Sk)|^2\,d\mu(\Sk)
\nonumber \\ & = &
\mu(X_{T,V} \bigtriangleup X_{Q^nx+T,V}) \nonumber \\ & = &
2[\mu(X_{T,V})-\mu(X_{T,V}\cap X_{Q^nx+T,V})]. \nonumber
\ee
The last equality uses the fact that $\mu$ is translation-invariant.
It follows that
\be \label{eq-conv2}
\mu(X_{T,V} \cap X_{Q^nx+T,V}) \to \mu(X_{T,V})
\ee
for each $\Tk$-tile $T$.

Now, combining the Lemma \ref{cyl-int-patch}, (\ref{eq-conv2}), and
Corollary \ref{cor-meas2}, we obtain that
$$
\freq(T\cup (Q^nx+T),\Tk) \to \freq(T,\Tk),\ \ \ \mbox{as}\ \ n\to\infty,
$$
for any tile $T\in \Tk$. In view of (\ref{eq-dens}) this implies
that $\lim_{n\to\infty}\dens(D_{Q^nx}) = 1$, and the proposition is proved. \qed

\medskip

\begin{defi}
A set $\Lam \subset \R^d$ is a Meyer set if $\Lam$ is Delone and
there is a finite set $F$ so that $\Lam - \Lam \subset \Lam + F$.
\end{defi}

As the next theorem shows,
under the additional assumption that $\Xi(\Tk)$ is a Meyer set,
the converse of Proposition~\ref{prop-pure} is also true.
Notice that we do not add
FLC to the assumption, since $\Xi(\Tk)$ being a Meyer set implies it.

\begin{theorem}\label{thm-pure}
Suppose that $\Tk$ is a repetitive fixed
point of a primitive substitution with expansive map $Q$
such that $\Xi(\Tk)$ is a Meyer set.
Then $(\Xt,\mu,\R^d)$ has pure discrete spectrum if and only if
\[
\lim_{n \to \infty} \dens(D_{Q^nx})=1, \ \ \ \mbox{for all} \ x \in \Xi(\Tk).
\]
\end{theorem}

This theorem extends \cite[Th.6.2]{soltil} to the case of $d \ge 3$.
Sufficiency was proved in Proposition~\ref{prop-pure}.
Necessity is not needed for our main results, so its proof (which is
a modification of the proof in \cite{soltil}) is given in the Appendix.

\medskip

\begin{remark} \label{rem-repet}
{\em In this subsection we have been assuming  that the fixed point
$\Tk$ of our
primitive substitution is {\em repetitive.} We note that this assumption
is not restrictive if we are interested in the measure-theoretic
dynamical system $(\Xt,\mu,\R^d)$. Indeed, if $\Tk$ is non-repetitive,
then the inclusion mapping $X_{\Ak,\om} \hookrightarrow \Xt$ induces an
isomorphism of measure-preserving systems (this follows from the proofs in
Subsection 6.1). All the tilings in
$X_{\Ak,\om}$ are repetitive; further, we can find a periodic point for
some power of $\om$ and work with it.}
\end{remark}

\subsection{Pure-pointedness on Delone multisets}
Transferring the results of the previous subsection over to Delone multisets, we obtain three quick
corollaries.

\begin{cor} \label{cor-erg}
If $\Lb$ is a primitive substitution Delone multiset with FLC such that
every $\Lb$-cluster is legal,
then the dynamical system $(X_{\Lbs}, \R^d)$ is uniquely ergodic.
\end{cor}

\noindent {\em Proof. } Apply Corollary \ref{leg-rep} and Lemma \ref{lem-conj}
to Theorem \ref{thm-erg}. \qed

\smallskip

\begin{cor} \label{cor-pure1}
Suppose that $\Lb=(\Lambda_i)_{i\le m}$ is a primitive substitution
Delone multiset with expansive map $Q$ such that
$\Xi(\Lb) := \cup (\Lam_i- \Lam_i)_{i \le m}$ is a Meyer set and
every $\Lb$-cluster is legal.
Then
the dynamical system $(X_{\Lbs},\mu,\R^d)$ has pure discrete
spectrum if and only if
\be
\lim_{n\to\infty} \dens(\Lb \cap (Q^nx+\Lb))=1,\ \ \
\mbox{for all}\ \  x \in \Xi(\Lb). \nonumber
\ee
\end{cor}

\noindent {\em Proof.} By Theorem \ref{legal-rep} $\Lb$ is representable by
a tiling $\Tk$. The legality of $\Lb$ implies that of $\Tk$ and
Lemma \ref{rep-legal} shows that $\Tk$ is repetitive.
The result is now a direct consequence of Theorem \ref{thm-pure}, in view of
the fact that $\Xi(\Tk) = \Xi(\Lb)$. \qed

\begin{cor} \label{cor-pure}
Suppose that $\Lb=(\Lambda_i)_{i\le m}$ is a primitive substitution Delone
multiset with expansive map $Q$ such that
$\bigcup_{i=1}^{m} \Lam_i$ lies in a lattice $L$ in $\R^d$ and every
$\Lb$-cluster is legal.
Then the dynamical system $(X_{\Lbs},\mu,\R^d)$ has pure discrete spectrum
if and only if
\be \label{eq-pure2}
\lim_{n\to\infty} \dens(\Lb \triangle (Q^nx+\Lb))=0,\ \ \
\mbox{for all}\ x \in L',
\ee
where $L'=L_1+\cdots+L_m$, and $L_i$ is the Abelian group generated by
$\Lam_i - \Lam_i$ for $ i \le m$.
\end{cor}

\noindent {\em Proof.} The Meyer set condition is obvious, since all the
sets $\Lam_i$ lie in the lattice $L$. We have $\Xi(\Lb) =
\cup (\Lam_i- \Lam_i)_{i \le m} \subset L'$,
so the necessity follows from Corollary~\ref{cor-pure1}. (This direction
will not be needed for our main results in the next section though.)
The condition (\ref{eq-pure2}) for $x\in \Lam_i-\Lam_i$ follows from
Corollary~\ref{cor-pure1}. In order to prove it for all $x\in L'$ we note
that $\Lb \triangle (y+z+\Lb) \subset (\Lb \triangle (y+\Lb))
\cup ((y+\Lb) \triangle (y+z+\Lb))$ for any $y, z \in \R^d$, hence
$$
\dens(\Lb \triangle (y+z+\Lb)) \le \dens(\Lb \triangle (y+\Lb))
+ \dens(\Lb \triangle (z+\Lb)),
$$
and the statement follows.
\qed

\section{Applications to lattice substitutions}

\subsection{Lattice substitutions}

In the sequel, $L$ will be a lattice in $\R^{d}$ and
the mappings of $\Phi$
will always be affine linear mappings of the form $x \mapsto Qx + a$,
where $Q \in \mbox{End}_{\Z}(L)$ is the {\em same} for all the maps and
$a \in L$.
Such maps are restrictions of uniquely determined affine linear mappings
on $\R^d$ and we will not distinguish them notationally.
A mapping $Q \in \mbox{End}_{\Z}(L)$ is called an {\em inflation} for $L$ if
$\det Q \neq 0$ and $\bigcap_{k=0}^{\infty} Q^k L = \{ 0 \}$.
So an expansive map $Q$ for $L$ is an inflation.

\begin{defi}
A {\em substitution system on $L$} with inflation $Q$ is a
pair $(\Lb,\Phi)$ consisting of
\begin{itemize}
\item a Delone multiset $\Lb = (\Lambda_i)_{i \le m}$,
for which each $\Lambda_i$ is a subset of $L$ and all $\Lambda_i$ are
mutually disjoint,
and
\item an $m \times m$ MFS $\Phi~$ on $L$, for which
\be
\Lambda_{i} = \bigcup_{j \le m} \bigcup_{f \in \Phi_{ij}} f(\Lambda_{j}),
\ \ \ i \le m\,,
 \label{8}
\ee
where the maps of $\Phi$ are affine linear mappings of the form
$x \mapsto Qx + a$, $a \in L$, and the unions in (\ref{8}) are disjoint.
\end{itemize}
\end{defi}

The substitution system $(\Lb,\Phi)$ is {\em primitive} if $S(\Phi)$ is
primitive.
For any affine linear mapping $f : x \mapsto Qx + b$ on $L$ we denote
the translational part,\,$b$,\,of $f$ by $t(f)$.

\medskip

\noindent
Let $(\Lb,\Phi)$ be a primitive substitution system on
$L$ with inflation $Q$. Then $(\Lb,\Phi^\ell)$ is a primitive
substitution system on $L$ with inflation $Q^\ell$, so without
loss of generality we can assume that $S(\Phi)$ is a positive matrix.
Let $L' := L_1 + L_2 + \cdots + L_m$, where $L_i := < \Lambda_i - \Lambda_i >$,
i.e.\ the Abelian group generated by $\Lambda_i - \Lambda_i$.

Since $Q$ is an inflation for $L$, $Q$ is an inflation for $L'$ also.
\noindent
In fact, \\
${\bigcup_{j \le m}\bigcup_{f \in \Phi_{ij}}(Q(\Lambda_j)+t(f))=\Lambda_i}$
implies
$\bigcup_{j \le m} Q(\Lambda_j - \Lambda_j) \subset (\Lambda_i - \Lambda_i)$  for any $i \leq m$, and thus
\be \label{propertyOfInflation}
Q(L_1 + \cdots + L_m) \subset L_i,\ \ \ \mbox{for any}\  i \leq m \,.
\ee
So $QL' \subset L'$, hence
$Q \in \mbox{End}_{\Z}(L')$. Note also
 $\bigcap_{k=0}^{\infty}
Q^k L' \subset \bigcap_{k=0}^{\infty} Q^k L = \{0\} \,.$

Note that $L/L'$ is finite, since $\Lb$ is a Delone multiset.
Let $q := |\det Q| = [L':QL'] > 1$.
We define the $Q$-adic completion
\[
\overline{L} = (\overline{L})_Q = \lim_{\leftarrow k} L/Q^kL'
= \lim_{\leftarrow k} (\dots \rightarrow L/Q^k L'
\rightarrow \dots \rightarrow L/QL'
\rightarrow L/L')
\]
of $L$ and
\[
\overline{L'} = (\overline{L'})_Q = \lim_{\leftarrow k} L'/Q^kL'
= \lim_{\leftarrow k} (\dots \rightarrow L'/Q^k L'
\rightarrow \dots \rightarrow L'/QL')
\]
of $L'$. Each of $\overline{L}$ and $\overline{L'}$ will be supplied with the
usual topology of a profinite group.
We can identify $\overline{L}$ and $(L/L') \times \overline{L'}$
as topological spaces.
In particular, the cosets $a + Q^k \overline{L'},
~a \in L, k = 1,2,\dots,$ form a basis of open sets of
$\overline{L}$ and each of these cosets is open and closed.
When we use the word coset in this paper, we mean either a coset of the form
$a + Q^k \overline{L'}$ in $\overline{L}$ or $a + Q^k L'$ in
$L$ according to the context.
An important observation is that any two cosets in $\overline{L}$
are either disjoint or one is contained in the other.
The same applies to cosets of $L$.

We let $\mu$ denote the Haar measure on $\overline{L}$,
normalized so that $\mu(\overline{L}) = 1$.
Thus for cosets, $\mu(a + Q^k \overline{L'}) = \frac{1}{|\det Q|^k
\cdot |L/L'|} = \frac{1}{q^k \cdot |L/L'|},~ k = 1,2,\dots$.
From $\bigcap_{k=0}^{\infty} Q^k L' = \{0\}$, we conclude that the mapping
$x \rightarrow \{ x ~\mbox{mod}~ Q^k L'\}_k$ embeds $L$
in $\overline{L}$.
We identify $L$ with its image in $\overline{L}$ and note that
$\overline{L}$ is then the closure of $L$.
With this identification, $L$ is a dense subgroup of
$\overline{L}$, so we have a unique extension of $\Phi$ to a MFS
on $\overline{L}$.
% in the obvious way.
Thus if $f \in \Phi_{ij}$ and $f: x \mapsto Qx + a$, this formula
defines a mapping on $\overline{L}$, to which we give
the same name.
Furthermore defining the compact subsets in $\overline{L}$
\[ W_i := \overline{\Lambda_i}, \ \ \ i \le m\,,
\]
and using the relations (\ref{8}) and the continuity of mappings, we have
\be \label{windowEquation}
  W_i = \bigcup_{j \le m} \bigcup_{f \in \Phi_{ij}} f(W_j),
\ \ \ i \le m \,.
\ee
We call $(\Wb,\Phi)$ the associated $Q$-adic system.

Suppose $L = \bigcup_{i \le m} \Lambda_i$.
For any $i \le m$, since $< \Lambda_i - \Lambda_i > \,\subset L'$,  we have
\[
\Lambda_i \subset x + L' \ \ \ \mbox{for any} \ x \in \Lambda_i \subset L.
\]
For $a \in L$, let
\be \label{def-a-class}
\Phi_{ij}[a]  &:=& \{ f \in (\Phi)_{ij} :~ Qy + t(f) \equiv a \mbox{~mod~} QL',\\
    \nonumber  &  &\mbox{where~} f : x \mapsto Qx+t(f),\
                     \Lambda_j \subset y + L' \} \\ \nonumber
              & =& \{ f \in (\Phi)_{ij} : f(\Lam_j) \subset a + QL'\}
\ee
Then
\[
\bigcup_{i,j \le m} \bigcup_{f \in \Phi_{ij}[a]} f(\Lambda_j) = a + QL'.
\]
Let $\Phi[a] := \cup_{i,j \le m} \Phi_{ij}[a]$.
This partitions $\Phi$ into congruence classes induced by $L/QL'$.

\begin{defi} \label{def-mod-coin}
Let $(\Lb,\Phi)$ be a primitive substitution system on $L$ with inflation
$Q$ and $L= \bigcup_{i \le m} \Lam_i$.
We say that $(\Lb,\Phi)$ admits {\em a modular coincidence
relative to $QL'$} if $\Phi[a]$ is contained entirely in one row of
$\Phi$ for some $a \in L$.
\end{defi}

It is easy to see that  $(\Lb,\Phi)$ admits a modular coincidence
relative to $QL'$ if and only if $(a+QL') \subset \Lambda_i$ for some
$i\le m$.

We define the notion of a model set (or cut and project set) \cite{RVM2}.
\begin{defi}
A {\em {cut and project scheme}} (CPS) consists of a collection of spaces and mappings as follows;
\be
\begin{array}{ccccc}
 \R^{d} & \stackrel{\pi_{1}}{\longleftarrow} & \R^{d} \times G & \stackrel{\pi_{2}}
{\longrightarrow} & G \\
 && \bigcup \\
 && \widetilde{L}
\end{array}
\ee
where $\R^{d}$ is a real Euclidean space, $G$ is some locally
compact Abelian group,
$ \widetilde{L} \subset {\R^{d}
\times G}$ is a lattice,  i.e.\
a discrete subgroup for which the quotient group
$(\R^{d} \times G) / \widetilde{L}$ is
compact, $\pi_{1}|_{ \widetilde{L}}$ is injective,
and $\pi_{2}(\widetilde{L})$ is dense in $G$.
\end{defi}
\begin{defi}
A {\em model set} in $\R^{d}$ is a subset of $\R^{d}$ which, up to
translation, is of the
 form $\Gam(V) = \{ \pi_{1}(x) : x \in \widetilde{L}, \pi_{2}(x) \in V\}$
for some cut and project scheme as above, where $V \subset G$ has non-empty
interior and compact closure (relatively compact).
The model set $\Gam(V)$ is {\em regular} if the boundary
$\partial V = \overline{V} \backslash \stackrel{\circ}{V}$ of $V$
is of (Haar) measure $0$.
\end{defi}
When we need to be more
precise we explicitly mention the cut and project scheme from which
a model set arises. This is quite important in some of the theorems
below. Model sets are always Delone subsets of $\R^d$.
We will also find it convenient to consider certain degenerate types of model
sets. A {\em weak} model set is a set in $\R^d$ of the form
$\Gam(V)$  where we assume only that $V$ is relatively compact, but
not that it has a non-empty interior. When $V$ has no interior,
$\Gam(V)$ is not necessarily relatively dense
in $\R^d$ but regularity still means that the boundary of $V$ is of measure $0$.

\subsection{Pure pointedness, modular coincidence, and model sets}

This subsection contains our main new result, namely, that
in the setting of lattice substitution systems, pure point diffraction
spectrum implies that there is a model set realization. Precise conditions
are given in Theorem~\ref{th-main}, which incorporates earlier results
and completes the circle of equivalences started in \cite{LM}.
The key new ingredient of the proof is
Theorem~\ref{th-modul}. Some of the arguments in this subsection are
similar to the corresponding parts of \cite{LM}.
However, there is an important distinction: here we have to do
everything modulo the sublattice $L'$. For instance, the notion of
modular coincidence is not the same as in \cite{LM}, $\overline{L}$ is
different, etc.

\smallskip

Consider $(\Lb,\Phi)$ a primitive substitution system on $L$ with an
expansive map $Q$.
Since $Q$ is expansive (see (\ref{def-expan})), for any bounded
subset $S$ of $\R^d$ containing $0$
as an interior point and $\lambda > 1$, there is large enough $k_0$
so that $Q^k S \supset \lambda S$ for $k > k_0$.

Recall that for any set $F\subset \R^d$ we have a cluster $F\cap \Lb=
(F\cap \Lambda_i)_{i\le m} \subset \Lb$.

Suppose $0 \in \Lb$.
Let $\Ib = (I_i)_{i \le m} := \Phi(0) = (\Phi_{ij}(0))_{i\le m}$
where $0 \in \Lambda_j$,  and let
$t(\Phi) := \{t(f) : f \in \Phi_{ij},~ f : x \mapsto Qx + t(f),
\ i, j \leq m \}$.
For a cluster $\Pb =(P_i)_{i\le m}$ we write
$\supp(\Pb)=\cup_{i\le m} P_i$.
Let $\{\beta_i \,:\, i = 1, \dots, d\}$ be a basis of $L'$.
Let $D_0$ be the parallelepiped in $\R^d$ of the form
\be
D_0 = \{ x_1 \beta_1 + \cdots + x_i \beta_i + \cdots + x_d \beta_d :
\, -1 \leq x_i < 1, 1 \leq i \leq d \,\}. \label{parallelepiped}
\ee
We can always find $p \in \Z_+$ so that $S(\Phi^p)$ is positive
and $Q^p D_0 \supset \lambda D_0$ for some $\lambda > 1$, and then
$a>0$ so that $D := a D_0 \supset (\supp(\Phi^p(0)) \cup t(\Phi^p))$.
Replacing $\Phi$ by $\Phi^p$, $QD_0$ by $Q^pD_0$ etc., we may, for the
purposes of Theorem \ref{th-modul}, assume at the outset that
$S(\Phi)$ is positive, $Q(D) \supset \lambda D$, and
$D \supset (\supp(\Ib) \cup t(\Phi))$ for some $\lambda > 1$ and convex $D$.

\medskip
\begin{lemma} \label{nice initial set}
Let $(\Lb,\Phi)$ be a primitive substitution system on $L$ with expansive
map $Q$ and $0 \in \Lb$.
Let $D$ be a convex set for which $Q(D)
 \supset \lambda D$ for some $\lambda > 1$
and $D \supset (\supp(\Ib) \cup t(\Phi))$.
Then there is $r > 0$ such that $\supp{(\Phi^n(\Ib))} \subset Q^n(rD)$
for all $n \in \Z_+$.
\end{lemma}

\noindent {\em Proof.}
Choose $p \in \Z_+$ so that $\lambda^p \geq p + 1$.
Since $\lambda D \subset Q(D)$, we have
$(p + 1)D \subset \lambda^p D \subset Q^{p}(D)$.
So
\be
(p+1)D \subset Q^p (D). \label{containmentOfD}
\ee
Note that  $Q(D) \supset D$ and
 $kD = \underbrace{D + \cdots + D}_{k}$ for any $k \in \Z_+$,
since $D$ is convex.\\
Since we have $\supp(\Ib) \subset D$, which means $\Ib \subset D\cap \Lb$,
for any $n \in \Z_+$
\be \label{eq-conta1}
\lefteqn{\supp{(\Phi^n(\Ib))}}\nonumber \\ \nonumber
& \subset & \supp{(\Phi^n(D\cap \Lb))}\\ \nonumber
& \subset & t(\Phi)+Q(\supp{(\Phi^{n-1}(D\cap\Lb))})  \\ \nonumber
& \subset & t(\Phi)+ \cdots+Q^{n-1}(t(\Phi)) +Q^n(D) \\
& \subset & D+\cdots +Q^{n-1}(D) +Q^n(D).
\ee
Since $Q(D)\supset \lam D\supset D$, we have $Q^i(D)\supset Q^j(D)$ for
$i>j$. Thus,
writing $n = lp + s$, where $1 \leq s \leq p$ and $l \in \Z_{\ge 0}$,
we obtain
\be
\sum_{i=0}^n Q^i(D)
& = & \sum_{i=0}^{n-lp} Q^i(D) + \sum_{k=1}^l Q^{n-kp} \lt(\sum_{i=1}^p
Q^i(D)\rt) \nonumber  \\ \nonumber
& \subset & (p+1) Q^{n-lp}(D) + \sum_{k=1}^l Q^{n-kp} (pQ^p(D)) \\ \nonumber
& = & Q^{n-lp}\lt( (p+1)D + p \sum_{k=1}^l Q^{(l-k+1)p}(D)\rt) \\ \nonumber
& \subset & Q^{n-lp}\lt( Q^p(D) + p\sum_{j=1}^l Q^{jp}(D)\rt) \ \ \
             {\mbox{from (\ref{containmentOfD}).}}\\ \nonumber
\ee
So
\be
\sum_{i=0}^n Q^i(D)
& \subset & Q^{n-lp} Q^p \lt( (p+1)D + p \sum_{j=1}^{l-1} Q^{jp}(D)\rt)
\nonumber \\ \nonumber
& \subset & Q^{n-lp} Q^p \lt( Q^p(D) + p \sum_{j=1}^{l-1} Q^{jp}(D)\rt)
\\ \nonumber
&  &        {\vdots}\\  \nonumber
& \subset & Q^{n-lp} \overbrace{Q^p \cdots Q^p}^l((p+1)D) = Q^n((p+1)D).
 \nonumber
\ee
In view of (\ref{eq-conta1}) this implies that $r:=p+1$ satisfies the
assertion of the lemma. \qed

\begin{theorem} \label{th-modul}
Let $(\Lb,\Phi)$ be a repetitive primitive substitution system on $L$
with expansive map $Q$ and $L = \bigcup_{i \le m} \Lambda_i$.
If $\dens(\Lb \,\triangle \,(Q^n{\alpha}+\Lb))
\stackrel{n \rightarrow \infty}{\longrightarrow} 0$\
for all $\alpha \in L'$,
then a modular coincidence relative to $Q^M L'$ occurs in $\Phi^M$
for some $M \in \Z_+$.
\end{theorem}

\noindent {\em Proof.}
Suppose that for all  $n \in \N$, $\Phi^n$ does not admit any modular
coincidence relative to $Q^n L'$.
We assume that $S(\Phi)$ is a positive matrix
without loss of the generality. Then the cluster
$\Ib=(I_i)_{i\le m} = (\Phi_{ij}(0))_{i \le m}$, with $0 \in \Lam_j$, has at
least one element
from each point set in $\Lb = (\Lambda_i)_{i \le m}$.
We claim that for all $n \in \Z_+$, $\supp(\Phi^n(\Ib))$ intersects
every coset $x+Q^nL'$ of $Q^nL'$ in $L$ non-trivially.
Indeed, $L=\bigcup_{i \le m} \Lam_i$ and each $\Lam_i \subset I_i+L'$ for
$i \le m$. So
\be
{L} & = & {\bigcup_{j \le m}\bigcup_{i \le m}(\Phi^n)_{ji}(I_i+L')}
\nonumber \\
  & = &  {\bigcup_{j \le m}\bigcup_{i \le m}((\Phi^n)_{ji}(I_i) +Q^nL')}
\nonumber \\
  & = & {\supp(\Phi^n(\Ib)) + Q^nL'.} \nonumber
\ee
So the claim follows.

To test modular coincidence we need only know about the translation parts of
$\Phi^n$ and to which coset of $L'$ each $\Lam_i$ belongs. (see Definition
\ref{def-mod-coin} and (\ref{def-a-class}).) Since $\Ib$ has at least one
element of each color type $i = 1, \dots, m$, modular coincidence
can be tested on $\Phi^n(\Ib)$:
if $(a+ Q^nL') \cap \supp(\Phi^n(\Ib)) \subset \Lambda_i$ for some $i$, then
it means all the mappings in $\Phi^n$ which contribute to produce points in
$a+Q^nL'$ lie on the $i$-row of $\Phi^n$.
So $\Phi^n$ has modular coincidence relative to $Q^nL'$.

Since $Q$ is an expansive map, there is a parallelepiped $D = a D_0$ for which
$Q(D) \supset \lambda D$ and $D \supset (\supp(\Ib) \cup t(\Phi))$ for
some $\lambda > 1$. Then by  Lemma \ref{nice initial set} there is
$r > 0$ such that $\supp{(\Phi^n(\Ib))}
\subset Q^n(rD)$ for all $n \in \Z_+$.

We have assumed that there is no modular coincidence. Thus, for any
$a\in L$ and any $n\ge 1$, there exists $i\le m$ and
\be \label{ur2}
x,y\in \supp(\Phi^n(\Ib))\cap (a+ Q^nL')\ \mbox{such that}\
x\in \Lambda_i,\  y\not\in \Lambda_i.
\ee
Let $\beta_{-j} = - \beta_j$ for $j=1,\ldots,d$.
We can write
$y-x$ as a non-negative integer linear combination of the vectors
$Q^n\beta_j,\ j\in \{\pm 1,\ldots, \pm d\}$. Now $Q^n(rD)$ is a
parallelepiped generated by $Q^n\beta_j$, $1 \le j \le d$,
 which contains $\supp(\Phi^n(\Ib))$. Then
there exists a path $x=x_1, x_2, \dots, x_s=y$ entirely in $Q^n(rD)$ whose
steps are each of the form $x_{l+1}-x_l=Q^n\beta_j, |j| \le d$ and we see that
there is a
$x' \in Q^n(rD) \cap (a+ Q^n L') \cap \Lambda_i \ \mbox{such that}\
x' + Q^n \beta_j \in  (Q^n(rD) \cap (a+ Q^n L')) \backslash \Lambda_i\ \
\mbox{for some}\ j, \ |j|\le d.$
It follows that
\be \label{ur3}
\sharp~ [ \cup_{|j|\le d} (\Lb \triangle
(Q^n{\beta_j}+\Lb)) \cap Q^n(rD)] \ge
| \det Q^n | \cdot |L/L'|,
\ee
since there are
$|\det Q^n| \cdot |L/L'|$ cosets of $Q^n L'$ in $L$ and
each of the cosets contributes at least one point to our count.

Furthermore for a parallelepiped $c +Q^n(rD)$ containing
$\supp(c+\Phi^n(\Ib))$,
for which $c + \Phi^n(\Ib)$ is a translate of $\Phi^n(\Ib)$,
the argument goes in the same way.
For $x,y$ in (\ref{ur2}) we have
$$
(x+c) \in \supp(c+\Phi^n(\Ib)) \cap \Lambda_i,\ \
(y+c) \in \supp(c+\Phi^n(\Ib)) \backslash \Lambda_i,\ \ y-x \in Q^n L'.
$$
As above, this implies that for some $x''\in
(c+ Q^n(rD)) \cap (c+a+ Q^nL') \cap \Lambda_i$ we
have $x'' + Q^n \beta_j \in  ((c+ Q^n(rD)) \cap (c+a+ Q^nL')) \backslash
\Lambda_i$ for some $j,\ |j|\le d$.
Thus, similarly to (\ref{ur3}),
\be \label{ur4}
\sharp~ [ \cup_{|j|\le d} (\Lb \triangle
(Q^n{\beta_j}+\Lb)) \cap (c+Q^n(rD))] \ge | \det Q^n | \cdot |L/L'|.
\ee

\noindent Let
$$
H_{0,R} = \{x\in L:\, x + (rD\cap \Lb) \subset B_R(0)\cap \Lb\},
$$
and let $\widetilde{H}_{0,R}$ be
a maximal set of $x\in H_{0,R}$ such that
$(x+rD)$ are mutually disjoint.
Since $Q$ is invertible, $Q^n(x+rD)$ and $Q^n(y+rD)$ are
disjoint if and only if $x+rD$ and $y+rD$ are disjoint.
So $\{Q^n(x+rD)\cap \Lb:\ x\in \widetilde{H}_{0,R}\}$
is a set of disjoint $\Lb$-clusters, which need not be translates of
each other. We claim
that each of the clusters $Q^n(x+rD)\cap \Lb$, for $x\in \widetilde{H}_{0,R}$,
contains a translate of $\Phi^n(\Ib)$. Indeed, if $x\in \widetilde{H}_{0,R}$,
then
$x + \Ib \subset B_R(0) \cap \Lb$ since $\Ib \subset rD\cap \Lb$.
It follows that
$\Phi^n(x+\Ib) \subset \Phi^n(\Lb)=\Lb$, $\Phi^n(x+\Ib) = Q^nx + \Phi^n(\Ib)$
since $x+\Ib$ is a translate of the cluster $\Ib$, and
\be
\supp{(\Phi^n(x+\Ib))} &=& Q^nx + \supp{(\Phi^n(\Ib))} \nonumber \\
&\subset& Q^nx+Q^n(rD) \nonumber \\
&=&Q^n(x+rD). \nonumber
\ee
 This claim,
together with (\ref{ur4}), yields
\be \label{ur5}
\lefteqn{\sharp~ [ \cup_{|j|\le d} (\Lb \triangle
(Q^n{\beta_j}+\Lb)) \cap Q^n(B_R(0))]} \nonumber \\ &  \ge &
\sharp~ \widetilde{H}_{0,R} \cdot | \det Q^n | \cdot |L/L'|.
\ee
Recall that
$$
\freq(rD\cap \Lb,\Lb) = \lim_{R\to \infty}
\frac{\sharp~ H_{0,R}}{\Vol(B_R)}\,,
$$
see Definition~\ref{def-ucf}. Let $e := \diam(rD)$.
The set $x+rD$, for $x\in L$, can intersect
at most $\sharp~ ((rD)^{+e}\cap L)$ translates $y+rD$, $y\in L$. Thus,
$$
\sharp~ \widetilde{H}_{0,R} \ge \frac{\sharp~ H_{0,R}}
{\sharp~ ((rD)^{+e}\cap L)}\,.
$$
Combining this with (\ref{ur5}) we obtain
$$
\frac{\sharp~ [ \cup_{|j|\le d} (\Lb \triangle
(Q^n{\beta_j}+\Lb)) \cap Q^n(B_R(0))]}{\Vol(Q^n B_R)} \ge
\frac{\sharp~ H_{0,R}\cdot |L/L'|}{\Vol(B_R)\cdot \sharp~ ((rD)^{e}\cap L)}\,.
$$
\noindent Letting $R\to \infty$, we conclude that
\[ \sum_{|j|\le d} \dens\,(\Lb \,\triangle \,
(Q^n \beta_j+\Lb)) \geq
\frac{\freq(rD\cap \Lb,\Lb)}{\sharp((rD)^{e}\cap L)}
\cdot |L/L'| \,.
\]
By assumption, $\textrm{dens}\,(\Lb \,\triangle\, (Q^n{\alpha}+\Lb))
\stackrel{n \rightarrow \infty}{\longrightarrow} 0$~ for all $\alpha \in L'$,
and in particular
\[
\sum_{|j|\le d} \dens(\Lb \triangle
(Q^{n}\beta_j+\Lb)) \stackrel{n \rightarrow \infty}
{\longrightarrow} 0\,.
\]
But $\freq(rD\cap \Lb,\Lb) >0$ by the repetitive property.
This is a contradiction.
\hfill $\square$

\medskip

\begin{theorem} \label{mainPFtheorem}
  Let $(\Lb,\Phi)$ be a primitive substitution system on $L$
with inflation $Q$.
Let $(\Wb,\Phi)$ be the corresponding associated $Q$-adic system.
Suppose $\overline{L} = \bigcup_{i \le m} W_{i}$.  Then
\begin{itemize}
\item[{\em (i)}] $S(\Phi^{r}) = (S(\Phi))^{r},\ \ \ \mbox{for all}
\ r \geq 1$;

\item[{\em (ii)}] $\mu(W_{i}) =
\frac{1}{q^{r}} \sum_{j \le m} (S(\Phi^{r}))_{ij} \mu(W_{j}),\ \ \ \mbox{for all}\;\;
i \le m,\ r \geq 1$;

\item[{\em  (iii)}] $\stackrel{\circ}{W_{i}} \neq \emptyset$ and
$ \mu (\partial W_{i}) = 0,\ \ \ \mbox{for all}\ i \le m$\,.
\end{itemize}
\end{theorem}

\noindent {\em Proof.}
For every measurable set $E \subset \overline{L}$ and any
$ f \in \Phi_{ij}$ where $f : x \mapsto Qx + a$,
$~ \mu(f(E)) = \mu (a+Q(E)) = \frac{1}{|\det Q|} \mu(E)$.
In particular, $\mu(f(W_{j})) = \frac{1}{q} ~w_{j},
~ \mbox{where}~ w_{j} := \mu (W_{j})~ \mbox{and}~  q=|\det Q| $. We obtain
\[
w_{i} \leq \sum_{j=1}^{m} \frac{1}{q^{r}}~ (S(\Phi^{r})_{ij})~w_{j},
\ \ \ \mbox{for any} \ r \geq 1
\]
from (\ref{windowEquation}).

Let $w = (w_i)_{i \le m}$.
Since $\overline{L} = \bigcup_{i \le m} W_{i}$,
the Baire category theorem assures
that for at least one $i$,
\be
\stackrel{\circ}{W_{i}} \neq \emptyset \label{11}
\ee
and then the primitivity gives this for all $i$.
So  $w > 0$ and
\be
 w \leq \frac{1}{q^{r}} S(\Phi^{r}) w \leq \frac{1}{q^{r}} S(\Phi)^{r} w,
~~ \mbox{for any}~ r \geq 1 \,.
\ee
Since $S(\Phi)^r$ is primitive and the PF eigenvalue of $S(\Phi)^r$ is
$~q^r = |\det Q|^r$ by \cite{lawa},  we have from \cite[Lemma 1]{LM} that
\be w = \frac{1}{q^{r}}S(\Phi^{r}) w =  \frac{1}{q^{r}}S(\Phi)^{r} w,
 ~~\mbox{for any}~ r \geq 1\, .
\ee
The positivity of $w$ together with  $S(\Phi^{r}) \leq S(\Phi)^{r}$ shows that
$S(\Phi^{r})= S(\Phi)^{r}$. This proves (i) and (ii).

Fix any $i \le m$, let $\stackrel{\circ}{W_{i}}$ contain
a basis open
set $a + Q^{r} \overline{L'}$ with some $r \in \Z_{\geq 0}$ by (\ref{11}).
 Since $(\Lb,\Phi^{r})$ is a substitution system,
$~ a+ Q^{r}\overline{L'} \subset \, \stackrel{\circ}{W_{i}} \, \subset
W_{i} = \bigcup_{j \le m} (\Phi^{r})_{ij} W_{j}$.
In particular, $(a + Q^{r}\overline{L'}) \cap g(W_{k}) \neq \emptyset\,$
for some
$k \le m$ and some $g \in (\Phi^{r})_{ik}$.
However $~g(z + \overline{L'}) = Q^{r}(z + \overline{L'})+ t(g)$, where
$\Lambda_k \subset z + L'$.
So $(a + Q^{r}\overline{L'}) \cap (t(g) + Q^r z + Q^{r}\overline{L'})
\neq \emptyset$.
This means $a + Q^{r}\overline{L'} = t(g) + Q^r z + Q^{r}\overline{L'}$. Thus
\be
g(W_{k}) \subset g(z + \overline{L'}) = a + Q^{r}\overline{L'} \subset
\stackrel{\circ}{W_{i}}\,.\label{14}
\ee
For all $f \in (\Phi^{r})_{ij},~ j \le m$, $f$ is clearly an open map,
so $\bigcup_{j \le m}(\Phi^{r})_{ij}(\stackrel{\circ}{W_{j}}) \subset \,
\stackrel{\circ}{W_{i}}\,.$ Thus
\begin{eqnarray}
 {\partial W_{i} = W_{i} ~\backslash \stackrel{\circ}{W_{i}}}
& = &{\left( \bigcup_{j \le m}(\Phi^{r})_{ij}(W_{j}) \right)
      ~\backslash \stackrel{\circ}{W_{i}}} \nonumber\\
& \subset &{\bigcup_{j \le m} \left( (\Phi^{r})_{ij}(W_{j})
      ~\backslash~ (\Phi^{r})_{ij}(\stackrel{\circ}{W_{j}}) \right)}\nonumber\\
& \subset &{\bigcup_{j \le m}(\Phi^{r})_{ij}(\partial W_{j})\,.} \label{15}
\end{eqnarray}
Note that due to (\ref{14}) at least one $g$ in $(\Phi^{r})_{ij}$
does not contribute to the relation (\ref{15}).

Let $v_{i} := \mu (\partial W_{i}),\, i \le m,$ and
$v := (v_i)_{i \le m}$. So $~v \leq
\frac{1}{q^{r}} S(\Phi^{r}) v$. Actually, by what we just said,
\be
0 \leq v \leq \frac{1}{q^{r}}S' v \leq \frac{1}{q^{r}}S(\Phi^{r})v
= \frac{1}{q^{r}}S(\Phi)^{r}v,\label{16}
\ee
where $S' \leq S(\Phi)^{r}, S' \neq S(\Phi)^{r}$.
Now applying \cite[Lemma 1]{LM} again we obtain
equality throughout (\ref{16}).
But by \cite[Lemma 2]{LM} the eigenvalues of $\frac{1}{q^{r}}S'$ are
strictly less in absolute value
than the PF eigenvalue of
$\frac{1}{q^{r}}S(\Phi)^{r}$, which is $1$.  This forces $v = 0$,
and hence $\mu (\partial W_{i}) = 0,
i \le m$. \hfill $\square$

\begin{theorem} \label{th-int-disj}
 Let $(\Lb,\Phi)$ be a primitive substitution system on $L$
with inflation $Q$ and $L = \bigcup_{i \le m} \Lambda_i$.
 Let $W_i = \overline{\Lambda_i}$
in $\overline{L}$ for any $i \le m$.
If there is a modular coincidence relative to $Q^M L'$ in $\Phi^M$,
then $\stackrel{\circ}{W_s} \cap \stackrel{\circ}{W_t} = \emptyset$,
 where $s, t \leq m$, $s \neq t$.
\end{theorem}

\noindent {\em Proof.}
By assumption, there is $\Lambda_i$ such that
\be
 a + Q^{M}L' = \bigcup_{j \le m} \bigcup_{f \in (\Phi^M)_{ij}[a]}
f(\Lambda_{j}) \subset \Lambda_i\, \ \mbox{for some} \ a \in L.
\ee
Assume $ (\Phi^M)_{ik}[a] \neq \emptyset$ for some $k$.
Pick $f \in (\Phi^M)_{ik}[a]$, where $f : x \mapsto Q^M x + t(f)$.
Then $Q^M y + t(f) \equiv a ~\mbox{mod}~ Q^M L'$, where $\Lambda_k \subset y + L'$,
 and $f(\Lambda_k) \subset a + Q^M L'$.

Take any $s \le m$ and suppose $\Lambda_s \subset z + L'$.
Pick $g \in \Phi_{ks} \neq \emptyset$, where $g : x \mapsto Qx + t(g)$.
Then
\begin{eqnarray*}
{f \circ g (\Lambda_s)} & = & {Q^M(Q(\Lambda_s) +t(g)) +t(f)}\\
                  & = & {Q^{M+1}(\Lambda_s)+Q^M(t(g)) + t(f)}\\
                  & \subset & { Q^{M+1}(z)+ Q^M(t(g)) +t(f)
                                + Q^{M+1} L'}.
\end{eqnarray*}
Let $c :=  Q^{M+1}(z)+ Q^M(t(g))+t(f)$.
So
\[ f \circ g (\Lambda_s) \subset c + Q^{M+1} L'.
\]
Let $p := f \circ g \in (\Phi^{M+1})_{is}$.
Since  $f \circ g (\Lambda_s) \subset a + Q^{M} L' \subset \Lambda_i$,
\[ (a + Q^M L') \cap (c + Q^{M+1} L') \neq \emptyset.
\]
Thus
\[ c + Q^{M+1} L' \subset a + Q^{M} L' \subset \Lambda_i.
\]
Let
\[ H_j := \{h \in (\Phi^{M+1})_{ij} : h(\Lambda_j) \subset c + Q^{M+1}L'\}
, \ \ \ \mbox{where}\ j \leq m.
\]
So
\be \label{unionEquation}
c + Q^{M+1}L' =
   \bigcup_{j \le m} \bigcup_{h \in H_{j}} h(\Lambda_{j})\, .
\ee
Note that for any $j \le m$ and $h \in H_j$,
  $Q^{M+1} x + t(h) \equiv c ~\mbox{mod~} Q^{M+1}L'$,
 where $ \Lambda_j \subset x + L'$.
So we can write (\ref{unionEquation}) more explicitly as follows;
\be \label{basicDecomposition}
c + Q^{M+1}L' & = & \bigcup_{\stackrel{j \le m}{h \in H_j}}
 \{c+Q^{M+1}\alpha_{h} + Q^{M+1}\Lambda_{j} : \nonumber \\
             &  & ~~~~~~~~~~c+Q^{M+1}\alpha_{h} = t(h),
                 \mbox{where}~ \alpha_{h} \in L~\}\,.
\ee
So
\[ L' = \bigcup_{j \le m} \bigcup_{h \in H_{j}} (\alpha_h+\Lambda_{j}).
\]
Note $\alpha_p+\Lambda_s \subset L'$.
Separating off $\Lambda_s$, we get
\be \label{Separating}
-\alpha_p + L' = \Lambda_s \cup (\bigcup_{j \le m} \bigcup_{h \in H'_{j}}
(- \alpha_p+\alpha_h+\Lambda_{j})\, ,
\ee
where $H'_j := H_j$ if $j \neq s$, and $H'_s := H_s \backslash \{p\}$.
Note that the decompositions of (\ref{Separating}) are disjoint.
But we also know that $\Lambda_s$ and $\bigcup \{\Lambda_j : \Lambda_j \subset -\alpha_p + L',
j \neq s \}$ are disjoint.
So it follows that
\be \label{containment}
\bigcup \{ \Lambda_j : \Lambda_j \subset -\alpha_p + L', j \neq s \} \subset
\bigcup_{j \le m}
\bigcup_{h \in H'_{j} }(- \alpha_{p}+\alpha_{h}+\Lambda_{j}).
\ee
Taking closures to both sides of (\ref{containment}),
\be \label{closureOfContain}
\bigcup \{ W_j : \Lambda_j \subset -\alpha_p + L', j \neq s \} \subset
\bigcup_{j \le m}
\bigcup_{h \in H'_{j} }(- \alpha_{p}+ \alpha_{h}+W_{j}).
\ee
On the other hand, if we apply Theorem \ref{mainPFtheorem}
 to $\Phi^{M+1}$ and look at
(\ref{basicDecomposition}), we get
 \[ \mu(c + Q^{M+1}\overline{L'}) = \sum_{j \le m}\sum_{h \in H_{j}}
\mu(Q^{M+1}(\alpha_{h}+W_{j}) + c)\, .
\]
Hence
\[ \mu(\overline{L'}) = \sum_{j \le m}\sum_{h \in H_{j}}
\mu (\alpha_{h}+W_{j})\, .
\]
So
\[ \mu(-\alpha_p + \overline{L'}) = \sum_{j \le m}\sum_{h \in H_{j}}
\mu (- \alpha_{p}+\alpha_{h}+W_{j})\, .
\]
Thus
\[ \mu(-\alpha_p + \overline{L'}) = \mu(W_{s}) +  \left( \sum_{j \le m}
\sum_{h \in H'_{j} } - \alpha_{p}+ \alpha_{h}+\mu(W_{j}) \right)
\]
which, after taking closures in  (\ref{Separating}),  gives us
\be \label{windowDisjoint1}
\mu \left( W_{s} \cap \left( \bigcup_{j \le m} \bigcup_{h \in H'_{j} }
(- \alpha_{p}+\alpha_{h}+W_{j}) \right) \right) = 0\,.
\ee
Finally from (\ref{closureOfContain}) and (\ref{windowDisjoint1})  we obtain
\[ \mu(W_{s} \cap (\bigcup \{W_j : \Lambda_j \subset -\alpha_p + L',
 j \neq s \})) = 0.
\]
It shows that  $\stackrel{\circ}{W_{s}} \cap \stackrel{\circ}{W_{j}} =
\emptyset ~$ for any $j$ with $\Lambda_j \subset -\alpha_p + L'$,
$ j \neq s $.
It is easy to see that $\stackrel{\circ}{W_{s}} \cap \stackrel{\circ}{W_{k}} =
\emptyset$, where $\Lambda_k \not\subset -\alpha_p + L'$.
Since $s$ is arbitrary in $\{1,\dots,m \}$,
$\stackrel{\circ}{W_{s}} \cap \stackrel{\circ}{W_{t}} = \emptyset$
for all $s, t \le m , s \neq t$. \hfill $\square$

\begin{lemma} \label{3lemmas}
Let $\Lambda_{i}, i \le m,$ be point sets of the lattice $L$ in $\R^{d}$.
Let $Q$ be an inflation of $L$, and identify $L$ and
its image in $\overline{L}$. Define
 $W_{i} := \overline{\Lambda_{i}}$ in $\overline{L}$ and
$\Gam_{i} := W_{i} \cap L$.
\begin{itemize}

\item[{\em (i)}]
 If $\Lambda_1, \dots ,\Lambda_m$ are disjoint and $\mu(\overline{\Gam_{i}
\backslash \Lambda_{i}}) = 0~ \mbox{for all}~ i \le m$, then
$\stackrel{\circ}{W_{i}} \cap \stackrel{\circ}{W_{j}} = \emptyset$
for all $ i \neq j$.

\item[{\em (ii)}] If $L = \bigcup_{i \le m} \Lambda_{i}$ and
 $\stackrel{\circ}{W_{i}} \cap \stackrel{\circ}{W_{j}} = \emptyset$
for all $ i \neq j$, where $ i,j \le m$,
then $\Gam_{i} \backslash \Lambda_{i} \subset \bigcup_{j \le m}
\partial W_{j}~~ \mbox{for all}~ i \le m $.

\item[{\em (iii)}] If $\mu(\partial W_{i}) = 0~$ for all $i \le m$ and
  $\Gam_{i} \backslash \Lambda_{i} \subset \bigcup_{j \le m} \partial W_{j}$,
then $\mu(\overline{\Gam_{i} \backslash \Lambda_{i}}) = 0$.
\end{itemize}
\end{lemma}
\noindent {\em Proof.}
(i) Suppose there are $ i, j \le m$
with $\stackrel{\circ}{W_{i}} \cap \stackrel{\circ}{W_{j}} \neq \emptyset$.
We can choose $a \in
{(\stackrel{\circ}{W_{i}} \cap \stackrel{\circ}{W_{j}}) \cap L}$,
since $L$ is dense in $\overline{L}$ and
$\stackrel{\circ}{W_{i}} \cap \stackrel{\circ}{W_{j}}$ is open.
Choose $k \in \Z_{+}$ so
that $a + Q^{k}\overline{L'} \subset \stackrel{\circ}{W_{i}}
\cap \stackrel{\circ}{W_{j}}$.
Note that $a +Q^{k}L' \subset \Gam_{i} \cap \Gam_{j}$. Then
\begin{eqnarray*}
{\bigcup_{i \le m} (\Gam_{i} \backslash \Lambda_{i})}
& \supseteq &{\left( (a + Q^{k}L')~ \backslash~ \Lambda_{i} \right)
\cup \left( (a + Q^{k}L')~ \backslash~ \Lambda_{j}
\right)}\nonumber\\
& \supseteq &{(a + Q^{k}L') ~\backslash~ (\Lambda_{i} \cap \Lambda_{j})}\nonumber\\
& = &{ a + Q^{k}L',~~\mbox{since the }~\Lambda_{i},~i \le m,~
\mbox{are disjoint}\,.}\nonumber
\end{eqnarray*}
So
\begin{eqnarray*}
 {\sum_{i \le m} \mu (\overline{\Gam_{i} \backslash \Lambda_{i}})}
& \geq & {\mu (\bigcup_{i \le m} (\overline{\Gam_{i}
\backslash \Lambda_{i}}))}\nonumber\\
 & \geq & {\mu (a + Q^{k}\overline{L'})}\nonumber\\
 & > &{0\,,}\nonumber
\end{eqnarray*}
contrary to assumption.

(ii) Assume $\stackrel{\circ}{W_{i}} \cap
\stackrel{\circ}{W_{j}} = \emptyset$  for all $i \neq j$.
For any $i \le m$,
\begin{eqnarray*}
 {(\Gam_{i} \backslash \Lambda_{i})}
& \subset & {( \bigcup_{j \neq i} \Lambda_{j}) \cap W_{i},
~~~\mbox{since}~ L = \bigcup_{i \le m} \Lambda_{i} }\\
 & \subset &{\bigcup_{j \neq i} (W_{j} \cap W_{i})}\\
& \subset & {\bigcup_{j \le m} \partial W_{j},
 ~~~\mbox{since}~ \stackrel{\circ}{W_{i}} \cap
 \stackrel{\circ}{W_{j}} = \emptyset~~ \mbox{for all}~ i \neq j \,.}
\end{eqnarray*}

(iii) Obvious. \hfill $\square$

\medskip

We consider the following cut and project scheme;

\be \label{cps2}
\begin{array}{ccccc}
\R^{d} & \stackrel{\pi_{1}}{\longleftarrow} &
\R^{d} \times \overline{L} &
\stackrel{\pi_{2}}{\longrightarrow} & \overline{L}\\
&& \bigcup \\
L & \longleftarrow & \widetilde{L} & \longrightarrow & L\\
t & \longleftarrow & (t,t) & \longrightarrow & t
\end{array}
\ee
where $\widetilde{L} := \{~(t,t)~ : ~t \in L \}
\subset \R^{d} \times \overline{L}$.

In fact, $\overline{L}$ is a compact abelian group and
$\widetilde{L} \subset \R^d \times \overline{L}$ is a lattice,
 i.e. a discrete subgroup for which the quotient group
$ (\R^d \times \overline{L}) / \widetilde{L}$ is compact.
Furthermore, $\pi_1 |_{\widetilde L}$ is injective and $\pi_2(\widetilde {L})$
is dense in $\overline{L}$.

\begin{theorem} \label{regModelSet}
Let $\Lambda_{i}, i \le m,$ be disjoint point sets of the lattice $L$ in
$\R^{d}$. Identify $L$ and its image in $\overline{L}$.
Let $W_{i} := \overline{\Lambda_{i}}$ in $\overline{L}$
and $\Gam_{i} := W_{i} \cap L$.
Suppose that $\mu(\partial W_{i}) = 0$ for all $i \le m$.
\begin{itemize}

\item[{\em (i)}]
 If $\Gam_{i} \backslash \Lambda_{i} \subset \bigcup_{j \le m} \partial W_{j}$
 then, relative to the CPS {\rm (\ref{cps2})}, $\Lambda_{i}$ is a regular
weak model set when
 $\stackrel{\circ}{W_{i}}$ is empty, and $\Lambda_{i}$ is a regular model set when
 $\stackrel{\circ}{W_{i}}$ is non-empty.

\item[{\em (ii)}]
 If  $L = \bigcup_{j \le m} \Lambda_{j}$ and each $\Lambda_i$ is a regular model set,
then  $\Gam_{i} \backslash \Lambda_{i} \subset \bigcup_{j \le m} \partial W_{j}$
for all $i \le m$.
\end{itemize}
\end{theorem}

\noindent {\em  Proof.}
 (i) Assume that $\Gam_{i} \backslash \Lambda_{i} \subset
\bigcup_{j \le m} \partial W_{j}$ for all $i \le m$.
Since $\mu(\partial W_{i}) = 0$ for all $i \le m$,
\be \label{measW}
\mu(W_{i}) = \mu(\stackrel{\circ}{W_{i}}) =
\mu(\stackrel{\circ}{W_{i}} \backslash \bigcup_{j \le m} \partial W_{j}).
\ee
Since $\Gam_{i} = W_{i} \cap L$, we have $\Lambda_{i} = V_{i} \cap L$
where $V_{i} := W_{i} \backslash (\Gam_{i} \backslash \Lambda_{i})$.
Now $V_{i} \supset\, \stackrel{\circ}{W_{i}} \backslash
\bigcup_{j \le m} \partial W_{j}$ by assumption.
 From
$\stackrel{\circ}{W_{i}} \backslash \bigcup_{j \le m} \partial W_{j}
\subset \, \stackrel{\circ}{V_{i}}
\subset V_{i} \subset \overline{V_{i}} = W_{i}$
and (\ref{measW}), $\mu(\overline{V_{i}} \backslash
\stackrel{\circ}{V_{i}}) = 0$.
So $\Lambda_i$ is regular.
If $\stackrel{\circ}{W_{i}} = \emptyset$, then
$\stackrel{\circ}{V_{i}} = \emptyset$ also.
Thus $\Lambda_i$ is a regular weak model set.
On the other hand, for any $i$ with $\stackrel{\circ}{W_{i}} \neq \emptyset$,
we know that $\stackrel{\circ}{V_{i}} \neq \emptyset$ and
$\overline{V_{i}}$ is compact. It follows that $\Lambda_{i} = V_{i} \cap L$
is a regular model set for the CPS (\ref{cps2}).

(ii) Suppose that  $\stackrel{\circ}{V_{i}} \neq \emptyset$,
$\mu(\overline{V_{i}} \backslash \stackrel{\circ}{V_{i}}) = 0$,
where $\Lambda_{i} = V_{i} \cap L$, and $L = \bigcup_{j \le m} \Lambda_{j}$.
Then from $\overline{\Gam_{i} \backslash \Lambda_{i}}
= \overline{(W_i \cap L) \backslash (V_i \cap L)}
\subset \overline{W_{i} \backslash V_{i}}
\subset W_{i} \backslash \stackrel{\circ}{V_{i}}
= \overline{V_{i}} \backslash \stackrel{\circ}{V_{i}}$, we have
$\mu(\overline{\Gam_{i} \backslash \Lambda_{i}}) = 0$ for all $i \le m$.
By Lemma \ref{3lemmas} (i) and (ii) ,
$\Gam_{i} \backslash \Lambda_{i} \subset \bigcup_{j \le m} \partial W_{j}$
for all $i \le m$.
\hfill $\square$

\begin{theorem}\label{martinThm} {\em (Schlottmann \cite{martin})} If
$\Gam \subset \R^d$ is a regular model set,
then $\Gam$ has pure point diffraction spectrum, i.e. the Fourier transform of
its volume averaged autocorrelation measure is a pure point measure.
\end{theorem}

This theorem was established for real internal spaces by \cite{Hof} and in
full generality, as stated here, in \cite{martin}. For a new simpler proof
of this result see \cite{BM}. \hfill $\square$

\medskip
Gathering all the results, we can state the following theorem, which
says, in particular, that $\Lb$ has pure point diffraction spectrum if
and only if each $\Lambda_i$ is a regular model set with respect to
a certain cut-and-project scheme, and that generalizes to lattice
substitutions in $\R^d$ Dekking's well-known criterion for pure point diffractivity.

\begin{theorem} \label{th-main}
Let $\Lb$ be a Delone multiset with expansive
map $Q$ such that $(\Lb,\Phi)$ is a primitive substitution system,
$L = \bigcup_{i \le m} \Lam_i$ for some lattice $L$ in $\R^d$, and every
$\Lb$-cluster is legal.
 Let $L' = L_1+\cdots+L_m$, where $L_i = <\Lam_i-\Lam_i>$. The following are
equivalent:

{\em (i)} $\Lb$ has pure point diffraction spectrum;

{\em (ii)} $\Lb$ has pure point dynamical spectrum;

{\em (iii)} $\dens(\Lb \triangle (Q^n \alpha+\Lb))
\stackrel{n \to \infty}{\longrightarrow} 0$
for all $\alpha \in L'$;

{\em (iv)}
 A modular coincidence relative to $Q^M L'$ occurs in $\Phi^M$ for some $M$;

{\em (v)} Each $\Lam_i$ is a regular model set for $i \le m$, relative to the
CPS (\ref{cps2}).
\end{theorem}

\noindent
{\em Proof.} It is easy to see that $\Lb$ has FLC, since $\Lb$ lies in a
lattice $L$ in $\R^d$, and that $\Lb$ is repetitive and has UCF, since
every $\Lb$-cluster is legal.
So the proof goes as follows:

(i) $\Leftrightarrow$ (ii) \cite[Th. 3.2]{LMS1}.

(ii) $\Rightarrow$ (iii) Corollary \ref{cor-pure}, using only the sufficiency,
which follows from Proposition~\ref{prop-pure}.

(iii) $\Rightarrow$ (iv) Theorem \ref{th-modul}.

(iv) $\Rightarrow$ (v) Theorem \ref{mainPFtheorem} and \ref{th-int-disj}, Lemma \ref{3lemmas}, and Theorem \ref{regModelSet}(i).

(v) $\Rightarrow$ (i) Theorem \ref{martinThm}.
\hfill $\square$

\begin{remark}
{\em If  $\Lb$ is  representable, but not repetitive then it still may be possible to use the criteria of
Theorem \ref{th-main} to check for pure point diffractivity,
using the argument in Remark \ref{rem-repet}.
One only has to find another $\Gb \in X_{\Lbs}$ which
does satisfy the conditions in Theorem \ref{th-main},
since $\Lb$ is pure point
diffractive if and only if $\Gb$ is pure point diffractive (see
\cite[Th. 3.2]{LMS1}).}

\end{remark}

\begin{example}
(Substitution Delone multiset with modular coincidence)\\
{\em Consider a substitution defined by $a \rightarrow abc$,
$b \rightarrow dcb$, $c \rightarrow cda$, and $d \rightarrow dab$.
We can consider a corresponding MFS $\Phi$ as follows;
\be
\Phi = \left( \begin{array}{cccc}
       \{3x\} & \emptyset & \{3x+2\} &\{3x+1\} \\
       \{3x+1\}& \{3x+2\} & \emptyset & \{3x+2\} \\
       \{3x+2\}& \{3x+1\} & \{3x\} & \emptyset \\
       \emptyset & \{3x\} & \{3x+1\} & \{3x\}
       \end{array} \right) \nonumber
\ee
A Delone multiset $\Lb = (\Lam_a,\Lam_b,\Lam_c,\Lam_d)$ generated from
($\{0\}$, $\{-1\}$, $\emptyset$, $\emptyset$) is fixed under $\Phi$.
On the real line, $\Lb$ looks like

\begin{picture}(100,50)(0,15)
 \put(5,40){\line(1,0){350}}
 \put(20,40){\multiput(20,0)(20,0){15}{\circle*{.20}}}
 \put(15,30){$\cdots$}
 \put(40,30){$c$}
 \put(60,30){$d$}
 \put(80,30){$a$}
 \put(100,30){$d$}
 \put(120,30){$c$}
 \put(140,30){$b$}
 \put(160,30){$a$}
 \put(180,30){$b$}
 \put(200,30){$c$}
 \put(220,30){$d$}
 \put(240,30){$c$}
 \put(260,30){$b$}
 \put(280,30){$c$}
 \put(300,30){$d$}
 \put(320,30){$a$}
 \put(340,30){$\cdots$}
\put(15,43){$\cdots$}
 \put(40,43){-6}
 \put(60,43){-5}
 \put(80,43){-4}
 \put(100,43){-3}
 \put(120,43){-2}
 \put(140,43){-1}
 \put(160,43){0}
 \put(180,43){1}
 \put(200,43){2}
 \put(220,43){3}
 \put(240,43){4}
 \put(260,43){5}
 \put(280,43){6}
 \put(300,43){7}
 \put(320,43){8}
 \put(340,43){$\cdots$}
\end{picture}

Note that $<\Lam_{\alpha}-\Lam_{\alpha}~ :~ \alpha \in \{a,b,c,d\}>  = 2\Z$
and that  $\Lam_a \subset 2\Z$, $\Lam_b \subset 1+2\Z$, $\Lam_c \subset 2\Z$,
and $\Lam_d \subset 1+2\Z$.
\begin{eqnarray*}
\Phi \left( \begin{array}{l}
       \Lam_a \\ \Lam_b \\ \Lam_c \\ \Lam_d
       \end{array}  \right)
 \subset
\left( \begin{array}{ccccccc}
       6\Z & \cup & \emptyset & \cup & (6\Z+2) & \cup & (6\Z+4) \\
       (6\Z+1) & \cup & (6\Z+5) & \cup & \emptyset & \cup & (6\Z+5) \\
       (6\Z+2) & \cup & (6\Z+4) & \cup & 6\Z & \cup & \emptyset \\
       \emptyset & \cup & (6\Z+3) & \cup & (6\Z+1) & \cup & (6\Z+3)
       \end{array}  \right).
\end{eqnarray*}
Since each of $\Phi[3]$ and $\Phi[5]$ lies in one row of $\Phi$, the modular
coincidence is confirmed. By Theorem \ref{th-main} the sets $\Lambda_1, \Lambda_2, \Lambda_3, \Lambda_4$ are
pure point diffractive and regular model sets. }
\end{example}

\subsection{ Concluding Remarks}

We are still left with some open questions around Theorem \ref{th-main}.
The most important is the requirement of legality, which we use as a
connecting link to employ the results for tilings. Is it sufficient that $\Lb$ is repetitive?
Example \ref{ex-nonlegal-rep} illustrates a case in which legality fails, but
repetitivity holds.
One can observe that each of the two point sets of which it is comprised is
pure point diffractive, even though the asymptotic condition of (iii) in
Theorem \ref{th-main} is not satisfied for some elements of $L'$.

The entire context of Theorem \ref{th-main} is that of an underlying lattice. Are there computable
conditions (like modular coincidence) for substitution Delone sets which link to pure point diffractivity when we move away from the lattice environment?

Thirdly, what is the true relationship of the two concepts of pure point diffractivity and regular
model sets? Assuming, say, a strictly ergodic substitution system satisfying the
 Meyer set condition (which is automatic for lattices), are these two concepts equivalent?

\section{Appendix}

\subsection{Unique ergodicity}

In this section we show that
if $\Tk$ is a fixed point of a primitive substitution,
then the tiling dynamical system $(X_{\Tk}, \R^d)$ is uniquely ergodic.
This is a bit more general than \cite[Th.3.1]{soltil} where it was also
assumed that $\Tk$ is repetitive. The proof goes through establishing
the existence of
{\em uniform patch frequencies} (UPF), the analog of UCF, see Definition
\ref{def-ucf}.
We present complete details here, in part because
the proof of UPF was omitted in \cite{soltil}, but we should note that in
a similar but slightly different setting the existence of UPF was
established in \cite[Prop.1]{GH} (see also \cite[p.182]{GH} for references to
earlier results of this kind).
%boris

Throughout this section, we fix $\Tk$  --- a tiling satisfying
$\om(\Tk)=\Tk$ for a primitive tile-substitution $\om$.

\smallskip

\begin{lemma} \label{lem-prag} {\em (see, e.g., \cite[1.6]{prag})}
Every prototile in a primitive tile-substitution
has the boundary of zero Lebesgue measure.
\end{lemma}

\begin{lemma} \label{lem-Van}
Let $\Tk$ be a fixed point of a substitution with expansive map $Q$.
Then for any
tile $T \in \Tk$ with $T=(A,i)$ for some $i \le m$,
$\{Q^n A\}_{n \ge 1}$ is a van Hove sequence.
\end{lemma}

\noindent
{\em Proof.} This is pretty straightforward from Lemma \ref{lem-prag}.
\qed

\medskip

The following is proved in \cite{soltil}. Note that it does not require
repetitivity. Below PF is an abbreviation for ``Perron-Frobenius.''

\begin{cor}
\label{cor-pf} Let $\om$ be a primitive tile-substitution with prototiles
$T_i=(A_i,i)$, for $i \le m$, and expansive map $Q$. Then
the PF eigenvalue of the substitution matrix $S$ is $|\det(Q)|$ and the vector
$(\Vol(A_i))_{i\le m}$ is a left PF eigenvector. Thus,
$$
\lim_{n\to\infty} |\det(Q)|^{-n} (S^n)_{ij} = r_i\Vol(A_j),
$$
where $(r_i)_{i\le m}$ is the right PF eigenvector of $S$ such that
$\sum_{i=1}^m r_i \Vol(A_i) = 1$.
\end{cor}

\medskip

\noindent {\sc Notation.}
For a patch $P$ and a bounded set $F\subset \R^d$ denote
$$
L_P(F) = \sharp \{g\in \R^d:\ g+P\subset \Tk,\ (g+\supp(P))\subset F\}
$$
and
$$
N_P(F) = \sharp \{g\in \R^d:\ g+P\subset \Tk,\ (g+\supp(P)) \cap F \ne \es\}.
$$

Let $V_{\mmin}$ and $V_{\mmax}$ be the minimal and maximal volumes of
$\Tk$-tiles respectively, let $\|Q\|$ be the operator norm, and let
$t_{\mmax}$ be the maximal diameter of $\Tk$-tiles.

\smallskip

\noindent
{\sc Inflated Tilings.} Given a tiling $\Tk$ and an expansive map $Q'$ on
$\R^d$ we let $$Q' \Tk = \{(Q'(\supp(T)), l(T)):\ T\in \Tk\}.$$
In other words, we blow up the tiles and retain their labels.
Usually we will take $Q' = Q^k$.
If $\Tk$ is a fixed point of $\omega$, then
the tilings $Q^k\Tk,\ k=1,2,\ldots,$ form
an hierarchical family of order $k$ {\em super-tilings} in the sense that every
tile of a higher order tiling can be decomposed into tiles of a lower order
tiling.

\begin{lemma} \label{lem-vspom1}
Let $\Tk$ be a fixed point of a primitive substitution.
Let $F\subset \R^d$ be an arbitrary bounded set, and let $\{F_n\}_{n\ge 1}$ be
a van Hove sequence in $\R^d$. Then for any $\Tk$-patch $P$ and any
$h\in \R^d$ we have

\vspace{1mm}

{\em (i)} $\ L_P(F) \le c_1\Vol(F)$, where $c_1$ depends only on $\Tk$;

\vspace{1mm}

{\em (ii)} If $P$ is a legal patch, then
there are $c_2, n_0 > 0$ (depending on $P$ and $\{F_n\}_{n\ge 1}$,
but not on $h$) so that
$\ L_P(h+F_n) \ge c_2\Vol(F_n)$, for all $n\ge n_0$;

\vspace{1mm}

{\em (iii)} If $P$ is a legal patch, then
$\ \lim_{n\to\infty} N_P(\partial (h+F_n)) /L_P(h+F_n) = 0$
uniformly in $h$.

\end{lemma}

\noindent {\em Proof.} (i) Select any tile from the patch $P$. Then distinct
$\Tk$-patches equivalent to $P$ will have distinct selected tiles. Therefore,
\[\Vol(F) \ge L_P(F) V_{\mmin},\]
so we can take $c_1 = V_{\mmin}^{-1}$.

(ii) If $P$ is a legal patch, then its translate occurs in some patch
$\om^k(T_i)$. Let $\ell\in \Nat$ be such that $\om^\ell(T_j)$ contains tiles
of all types, for all $j\le m$. (This exists by the
primitivity of the substitution.) Then every patch $\om^{k+\ell}(T),\
T\in \Tk$, contains a translate of $P$. We consider the super-tiling
$Q^{k+\ell}\Tk$. It follows that for any set $F$,\ \
$L_P(F)$ is at least the number
of $Q^{k+\ell}\Tk$-tiles whose supports are contained in $F$.
Therefore, for $r=\|Q\|^{k+\ell}\cdot t_{\mmax}$,
\be
\lefteqn{L_P(h+F_n)\cdot |\det(Q)|^{k+\ell} V_{\mmax}} \nonumber \\ \nonumber
& \ge & \Vol(h+F_n^{-r}) =
\Vol(F_n^{-r}) \ge \Vol(F_n) - \Vol((\bd F_n)^{+r}).
\ee
This implies the desired statement in view of (\ref{Hove}).

(iii) Let $t = \diam(P)$.
$$ \frac{N_P(h+\bd F_n)}{L_P(h+F_n)} \le
\frac{L_P(h+{\bd F_n}^{+t})}{L_P(h+F_n)} \le
\frac{c_1 \Vol(h+{\bd F_n}^{+t})}{c_2 \Vol(h+F_n)} \rightarrow 0,\,
\mbox{from (i) and (ii)}.
$$
%(iv) This is immediate from (iii), (ii), and (i).
\qed

\medskip

\begin{lemma} \label{lem-freq}
Let $\Tk$ be a fixed point of a primitive substitution with expansive map $Q$.
Let $P$ be a $\Tk$-patch. Then
$$
c_P := \lim_{n\to\infty} \frac{L_P(Q^n A)}{\Vol(Q^n A)}
$$
exists uniformly in A, a support of a $\Tk$-tile.
\end{lemma}

\noindent {\em Proof.} If $P$ is non-legal, then $L_P(Q^n A)=0$ for every
$n$ and every tile support $A$, so $c_P=0$.

Assume now $P$ is a legal patch.
Let $\{T_1,\ldots,T_m\}$ be representatives of all $\Tk$-tile types,
having supports $A_i$, $i \le m$.
Fix $\e>0$. By Lemma~\ref{lem-vspom1}(iii) and Lemma~\ref{lem-Van},
we can find $k_0\in \Nat$ so that
for any $k \ge k_0$ and any tile support $A$ on $\Tk$,
\be \label{eq-fre1}
N_P(\bd Q^k A) \le \e L_P(Q^k A)\,.
\ee
Choose a tile $T \in \Tk$. Then $T=(A,j)$ for some $j \le m$, where
 $A= \supp(T)$.
Consider the subdivision of $Q^n A=\supp(Q^n T)$, $n> k >k_0$,
into the tiles
of $Q^k \Tk$. By the definition of the substitution matrix $S$, there are
$(S^{n-k})_{ij}$ tiles equivalent to $Q^k T_i$ in the subdivision of
$Q^n T$. Therefore, in view of (\ref{eq-fre1}),
\be \label{eq-fre2}
\sum_{i=1}^m L_P(Q^k A_i) (S^{n-k})_{ij} \le L_P(Q^n A) \le (1+\e)
\sum_{i=1}^m L_P(Q^k A_i) (S^{n-k})_{ij}.
\ee
By Corollary~\ref{cor-pf},
$$
\lim_{n\to\infty} \frac{(S^{n-k})_{il}}{\Vol(Q^n A_l)} = r_i|\det(Q)|^{-k},
\ \ \mbox{for any} \ \ l \le m.
$$
Thus, dividing (\ref{eq-fre2}) by $\Vol(Q^n A)$ and letting $n\to\infty$
we obtain
$$
\limsup_{n\to\infty} \frac{L_P(Q^n A)}{\Vol(Q^n A)} -
\liminf_{n\to\infty} \frac{L_P(Q^n A)}{\Vol(Q^n A)} \le \e |\det(Q)|^{-k}
\sum_{i=1}^m r_i L_P(Q^k A_i).
$$
By Lemma~\ref{lem-vspom1}(i), the right-hand side does not exceed
$c_1 \e\sum_{i=1}^m r_i \Vol(A_i) = c_1\e$. Since $\e>0$ is arbitrary,
this proves the existence of the limit.
Moreover,
$$
\sum_{i=1}^m L_P(Q^k A_i) r_i |\det(Q)|^{-k} \le c_P \le
(1+\e) \sum_{i=1}^m L_P(Q^k A_i) r_i |\det(Q)|^{-k}
$$
and it does not depend on the choice of the tile $T$ on $\Tk$.
Since $L_P(Q^k A_i)>0$ for $k$ sufficiently large,
we have $c_P>0$.
 \qed

\begin{lemma} \label{lemma-freq}
Let $\Tk$ be a fixed point of a primitive substitution.
For any $\Tk$-patch $P$ and for any van Hove sequence $\{F_n\}_{n \ge 1}$,
there exists
\be \label{eq-fri}
\freq(P,\Tk) := \lim_{n\to\infty} \frac{L_P(h+F_n)}{\Vol(F_n)} = c_P,
\ee
uniformly in $h\in \R^d$.
\end{lemma}

\noindent {\em Proof.}
Consider the decomposition of the space $\R^d$ into the tiles of $Q^k \Tk$ for $k$ large.
Then $L_P(h+F_n)$
is roughly the sum of $L_P(Q^k A)$ where $Q^k A$ ranges over the
supports of those $Q^k\Tk$-tiles which intersect $h+F_n$. For large $n$ the
``boundary effects'' from $\bd F_n$ become small by
the definition of van Hove sequence.
Note also that the ``boundary effects'' from
$\bd Q^k A$ become small.
%To be careful, one first needs to choose $k$, then $n$.
This is the idea; now let us give the details.

Let $G_{k,n}=\{A :\ Q^k A \cap (h+F_n) \neq \emptyset,\
A=\supp(T) \ \mbox{for}\ T\in \Tk\}$ and $H_{k,n} =\{A :\ Q^k A \subset
(h+F_n),\ A=\supp(T) \ \mbox{for}\ T\in \Tk\}$ for $k, n \ge 1$. We have
\be \label{esta1}
\sum_{A\in H_{k,n}} L_P(Q^k A) \le L_P(h+F_n) \le
\sum_{A\in G_{k,n}} [L_P(Q^k A) + N_P(\bd(Q^k A))].
\ee
Fix $\eps>0$. Using Lemma~\ref{lem-freq} and Lemma~\ref{lem-vspom1}(iii) choose
$k$ so that for any tile support $A$,
\be \label{esta2}
|L_P(Q^k A)/\Vol(Q^k A) - c_P| < \eps\ \ \mbox{and}\ \
N_P(\partial(Q^k A)) < \eps L_P(Q^k A).
\ee
Combining (\ref{esta1}) and (\ref{esta2}) we obtain
$$
(c_P-\eps) \sum_{A\in H_{k,n}} \Vol(Q^k A) \le L_P(h+F_n) \le
(1+\eps)(c_P+\eps) \sum_{A\in G_{k,n}} \Vol(Q^k A).
$$
Let $t_k:= \max\{\diam(Q^k A):\ A=\supp(T), T\in \Tk\}$. Observe that
$$
\sum_{A\in H_{k,n}} \Vol(Q^k A) \ge \Vol(F_n^{-t_k})\ \ \mbox{and}\ \
\sum_{A\in G_{k,n}} \Vol(Q^k A) \le \Vol(F_n^{+t_k}).
$$
Since $k$ is fixed,
by the van Hove property we have for $n$ sufficiently large:
$$
\Vol(F_n^{-t_k}) \ge (1-\eps) \Vol(F_n)\ \ \mbox{and}\ \
\Vol(F_n^{+t_k}) \le (1+\eps) \Vol(F_n).
$$
Combining everything, we obtain for $n$ sufficiently large:
$$
(c_P-\eps)(1-\eps) \le \frac{L_P(h+F_n)}{\Vol(F_n)} \le (c_P+\eps)(1+\eps)^2.
$$
Since $\eps$ was arbitrary, this implies (\ref{eq-fri}) as desired. \qed

\subsection{Overlap coincidence}

Here we prove the necessity in Theorem~\ref{thm-pure}.
We use the notion of overlaps and the
subdivision graph of overlaps from \cite[p.721]{soltil}, with some
modifications.

\begin{defi} \label{def-overlap} Let $\Tk$ be a tiling. A triple $(T,y,S)$,
with $T,S \in \Tk$ and $y \in \Xi(\Tk)$, is called an {\em overlap}
if the intersection
$\supp(y+T)\cap \supp(S)$ has non-empty interior.
We say that two overlaps $(T,y,S)$ and $(T',y',S')$ are {\em equivalent} if
for some $g\in \R^d$ we have $y+T = g+y'+T',\ S = g+S'$. Denote by
$[(T,y,S)]$ the equivalence class of an overlap. An overlap $(T,y,S)$ is a
{\em coincidence} if $y+T = S$. The {\em support} of an overlap
$(T,y,S)$ is $\supp(T,y,S) = \supp(y+T)\cap \supp(S)$.
\end{defi}

\begin{lemma} \label{lem-finite} Let $\Tk$ be a tiling
such that $\Xi(\Tk)$ is a Meyer set. Then the number of equivalence classes
of overlaps for $\Tk$ is finite.
\end{lemma}

\noindent {\em Proof.}
Let $T_i,\,i\le m$, be the representatives of all tile types for $\Tk$.
Let $\Lambda_i$ be the Delone set such that $\Lambda_i+T_i$ is the collection
 of all tiles of type $i$. Thus we have $\Tk = \bigcup_{i\le m}
(\Lambda_i+T_i)$.
Let $(T,y,S)$ be an overlap.
We have $T=u_i+T_i$ and $S = u_j+T_j$ for some
$i,j\le m$ (possibly equal) and some $u_i\in \Lambda_i,\,u_j\in \Lambda_j$. The
equivalence class of the overlap is completely determined by $i,j$, and
the vector $u_i + y- u_j$. Since
the interiors of the supports of $y+T$ and $S$ must intersect, we have
\be \label{eq-bound}
|u_i + y- u_j| \le C,
\ee
where $C= 2 \mmax \{\diam(T) : T \in \Tk\}$. Note that $u_i,
y, u_j\in \Xi$. By the definition of $\Xi:=\Xi(\Tk)$ we have
$\Xi = -\Xi$. By the definition of Meyer set, $\Xi-\Xi \subset \Xi+F$
for some finite set $F$. This implies $\Xi+\Xi-\Xi \subset \Xi + (F+F)$,
which is a discrete set, so there are finitely many possible vectors
$u_i + y- u_j$ in (\ref{eq-bound}). This proves the lemma. \qed

\medskip

Next we define the {\em subdivision graph $\Gk_{\Ok}(\Tk)$ for overlaps}.
Its vertices are equivalence classes of overlaps. Let $\Ok=(T,y,S)$ be
an overlap.  We will specify directed edges
leading from the equivalence class $[\Ok]$. Recall that we have
the tile-substitution $\om$, see Definition~\ref{def-subst}. Then
$\om(y+T) = Qy+\om(T)$ is a patch of $Qy+\Tk$, and $\om(S)$ is a $\Tk$-patch,
and moreover,
$$
\supp(Qy+\om(T)) \cap \supp(\om(S)) = Q(\supp(T,y,S)).
$$
For each pair of tiles $T'\in \om(T)$ and $S'\in \om(S)$ such that
$\Ok':= (T',Qy,S')$ is an overlap, we draw a directed edge from $[\Ok]$ to
$[\Ok']$.

\medskip

\begin{lemma} \label{lem-geom}
Let $\Tk$ be a repetitive fixed point of a primitive substitution with
expansive map $Q$ such that $\Xi(\Tk)$ is a Meyer set.
Let $x\in \Xi(\Tk)$. The following are equivalent:

 {\em (i)} $\lim_{n\to\infty} \dens(D_{Q^nx})=1$;

 {\em (ii)} $1-\dens(D_{Q^n x}) \le Cr^n$, $n\ge 1$,
for some $C > 0$ and $r\in (0,1)$;

 {\em (iii)} From each vertex of the graph $\Gk_{\Ok}(\Tk)$ there is a path
leading to a coincidence.
\end{lemma}

\noindent {\em Proof.}
We have for each $n\ge 0$:
$$
\R^d = \bigcup_{T\in \Tk} \bigcup_{S\in \Tk} \supp(T,Q^nx,S),
$$
where the support is considered empty if $(T,Q^nx,S)$ is not an overlap.
Notice that $\supp(D_{Q^nx})$ is exactly the union of supports of coincidences
in this formula. It is clear that all edges from the overlaps-coincidences
lead to other coincidences. Thus,
$$
Q(\supp(D_{Q^nx}))\subset \supp(D_{Q^{n+1}x}),
$$
and
\be \label{inc-dens}
\dens(D_{Q^nx}) \le \dens(D_{Q^{n+1}x}).
\ee
Now, if (iii) holds, then there exists $\ell\in\Nat$ such that
for each overlap $\Ok$ the ``inflation''
$Q^\ell(\supp(\Ok))$ contains a coincidence.
The volume of the support of a coincidence is at least $V_{\mmin}$. Thus
$$
1-\dens(D_{Q^{n+\ell}x}) \le \left(1 -
\frac{V_{\mmin}}{V_{\mmax}|\det(Q)|^\ell} \right)\left(1-\dens(D_{Q^nx})
\right).
$$
For any $n \ge 0$, $n = k\ell+s$ for some $k \in \N$ and $0 \le s < \ell$.
So
\be
1 - \dens(D_{Q^nx}) & = & 1 - \dens(D_{Q^{k\ell+s}x}) \nonumber\\
& \le & {b}^k (1 - \dens(D_{Q^s x})),\ \mbox{where}\
b = 1 - \frac{V_{\mmin}}{V_{\mmax}|\det(Q)|^\ell} \nonumber\\
& = & ({b}^{1/\ell})^{k\ell+s} \frac{(1-\dens(D_{Q^sx}))}{{b}^{s/\ell}}
\nonumber\\
& \le & r^n C, \ \mbox{for some}\ r \in (0,1) \ \mbox{and}\ C > 0.
\ee
Then (ii) follows.

It is straightforward that (ii) implies (i).

It remains to prove that (i) implies (iii).
Suppose to the contrary,
that there is an overlap $\Ok$ from which there is no path
to a coincidence. Then $Q^n(\supp(\Ok))\subset \R^d \setminus
\supp(D_{Q^nx})$ for all $n$. By the repetitivity of $\Tk$, the overlaps
equivalent to $\Ok$ occur relatively dense in $\R^d$. Therefore,
$$
1-\dens(D_{Q^nx}) \ge \dens(Q^n(\supp(\Ok))=\dens(\supp(\Ok)) >0,
$$
which contradicts (i). This completes the proof of the lemma.
\qed

\medskip

\noindent {\em Proof of Theorem~\ref{thm-pure}.}
The assumption gives us through Lemma \ref{lem-geom}
that there exists a basis ${\mathcal B}$ for $\R^d$ such that
for all $x \in {\mathcal B}$,
\be \label{sum-density}
\sum_{n=0}^{\infty}(1-\dens(D_{Q^nx})) < \infty.
\ee
Then \cite[Th.6.1]{soltil} implies that the dynamical system
$(\Xt,\mu,\R^d)$ has pure discrete spectrum. \qed

\end{document}